\documentclass[a4paper,11pt, reqno]{amsart}
\usepackage{latexsym}
\usepackage{amsbsy}
\usepackage{amsfonts}
\usepackage{amsmath}
\usepackage{geometry}
\usepackage[T1]{fontenc}
\usepackage{graphics,graphicx}
\usepackage [dvips]{epsfig}

\def\marginpar#1{\ignorespaces}

\newtheorem{thm}{Theorem}
\newtheorem{prop}[thm]{Proposition}

\newtheorem{definition}[thm]{Definition}

\theoremstyle{definition}

\newtheorem{remark}[thm]{Remark}

\numberwithin{equation}{section} \numberwithin{thm}{section}

\newdimen\AAdi%
\newbox\AAbo%
\def\AAk#1#2{\setbox\AAbo=\hbox{#2}\AAdi=\wd\AAbo\kern#1\AAdi{}}%

\def\eqref#1{(\ref{#1})}
\def\eqlabel#1{\def\@currentlabel{#1}}

\def\formula#1{\def\@tempa{#1}\let\@tempb\theequation\def\theequation{%
\hbox{#1}}\def\@currentlabel{(\theequation)}$$}
\def\endformula{\leqno\hbox{(\@tempa)}$$\@ignoretrue\let\theequation\@tempb}

\def\given{\hskip5\p@\relax\vrule\@width.4\p@\hskip5\p@\relax}

\newcommand{\open}[1]{%
\par\normalfont\topsep6\p@\@plus6\p@\trivlist\item[\hskip\labelsep\itshape#1%
\@addpunct{.}]\ignorespaces}

\DeclareRobustCommand{\close}[1]{%
  \ifmmode 
  \else \leavevmode\unskip\penalty9999 \hbox{}\nobreak\hfill
  \fi
  \quad\hbox{$#1$}}

\newlength{\toskip}

\def\BB{\mathcal B}
\def\CC{\mathcal C}

\def\FF{\mathcal F}

\def\NN{\mathcal N}

\def\SS{\mathcal S}

\def\vep{\varepsilon}
\def\<{\langle}
\def\>{\rangle}

\def \R {{\mathbb R}}

\def \T {{\mathbb T}}
\def \P {{\mathbb P}}
\def \E {{\mathbb E}}
\def \G {{\mathbb G}}
\def \H {{\mathbb H}}
\def \N {{\mathbb N}}

\makeatother

\begin{document}

\title[Deviation inequalities for bifurcating Markov chains on Galton-Watson tree]
{Deviation inequalities for bifurcating Markov chains on
Galton-Watson tree}

\author[S.V. Bitseki Penda]{{S.Val\`ere} Bitseki Penda }
\address{{\bf {Val\`ere} BITSEKI PENDA}\\ Laboratoire de Math\'ematiques,
CNRS UMR 6620, Universit\'e Blaise Pascal, avenue des Landais 63177 Aubi\`ere.}
\email{Valere.Bitsekipenda@math.univ-bpclermont.fr}

\begin{abstract}
We provide deviation inequalities for properly normalized sums of
bifurcating Markov chains on Galton-Watson tree. These processes are
extension of bifurcating Markov chains (which was introduced by
Guyon to detect cellular aging from cell lineage) in case the index
set is a binary Galton-Watson process. As application, we derive
deviation inequalities for the least-squares estimator of
autoregressive parameters of bifurcating autoregressive processes
with missing data. These processes allow, in case of cell division,
to take into account the cell's death. The results are obtained
under an uniform geometric ergodicity assumption of an embedded
Markov chain.
\end{abstract}

\maketitle \textit{Key words: Bifurcating Markov chains,
Galton-Watson processes, ergodicity, deviation inequalities, first
order bifurcating autoregressive process with missing data, cellular
aging}.

\vspace{2pt}


\section{Introduction}\label{gw_intro}
Bifurcating Markov chains (BMC) on Galton-Watson (GW) tree are an
extension of BMC to GW tree data. They were introduced by Delmas and
Marsalle \cite{DelMar} in order to take into account the death of
individuals in the Escherichia coli's (E.coli) reproduction model.
E.coli is a rod-shaped bacterium which reproduces by dividing in the
middle, thus producing two cells. One which has the new pole of the
mother and that we call new pole progeny cell, and the other which
has the old pole of the mother and that we call old pole progeny
cell. In fact, each daughter cell has two poles. One which is new
(new pole) and the other which already existed (old pole). The age
of a cell is given by the age of its old pole (i.e the number of
generations in the past of the cell before the old pole was
produced).

Guyon \& Al \cite{Gu&Al}  proposed the following linear Gaussian
model to describe the evolution of the growth rate of the population
of cells derived from an initial individual:
\begin{equation}\label{gw_bar11}
\begin{array}{ll}
\mathcal{L}(X_{1})=\nu, \qquad{\rm and}\qquad  \forall n\geq 1,
\quad \left\{\begin{array}{ll} X_{2n}=\alpha_{0}X_{n} + \beta_{0} +
\varepsilon_{2n} \\ \\ X_{2n+1} = \alpha_{1}X_{n} + \beta_{1} +
\varepsilon_{2n+1},
\end{array} \right.
\end{array}
\end{equation}
where $X_{n}$ is the growth rate of individual $n$, $n$ is the
mother of $2n$ (the new pole progeny cell) and $2n+1$ (the old pole
progeny cell), $\nu$ is a distribution probability on $\mathbb{R}$,
$\alpha_{0}, \alpha_{1}\in (-1,1)$; $\beta_{0}, \beta_{1}\in
\mathbb{R}$ and $\big((\varepsilon_{2n}, \varepsilon_{2n+1}), n\geq
1\big)$ forms a sequence of i.i.d bivariate random variables with
law $\mathcal{N}_{2}(0,\Gamma)$, where
\begin{equation}\label{cov_mat}
\Gamma=\sigma^{2}\begin{pmatrix} 1 & \rho \\ \rho & 1
\end{pmatrix}, \quad \sigma^{2}>0, \quad \rho\in(-1,1).
\end{equation}
The processes $(X_{n})$ defined by (\ref{gw_bar11}) are typical
example of BMC which are called the first order bifurcating
autoregressive processes (BAR(1)). The BAR(1)  processes are an
adaptation of autoregressive processes, when the data have a binary
tree structure (see Figure 1). They were first introduced by Cowan
and Staudte \cite{CS86} for cell lineage data where each individual
in one generation gives rise to two offspring in the next
generation.

In \cite{Guyon}, Guyon, using the theory of BMC, gave laws of large
numbers and central limit theorem for the least-squares estimator
$\widehat{\theta}^{r}=(\widehat{\alpha}_{0}^{r},
\widehat{\beta}_{0}^{r}, \widehat{\alpha}_{1}^{r},
\widehat{\beta}_{1}^{r})$ of the 4-dimensional parameter
$\theta=(\alpha_{0},\beta_{0},\alpha_{1},\beta_{1})$. He has also
built some statistical tests which allow to test if the model is
symmetric or not, and if the new pole and the old pole populations
are even distinct in mean. This allowed him to conclude a
statistical evidence in aging in E. Coli. Let us also mention
\cite{BerSapGég}, where Bercu \& Al. using the martingale approach
give asymptotic analysis of the least squares estimator of the
unknown parameters of a general asymmetric $p$th-order BAR
processes.

However, in the BMC model presented by Guyon, cells are assumed to
never die (a death corresponds to no more division). To take into
account cells's death, Delmas and Marsalle \cite{DelMar}, instead of
a regular binary tree, used a binary GW tree to label cells. In the
sequel, we will introduce the model which allowed them to study the
behavior of the growth rate of cells, taking into account their
possible death.

\subsection{The model}\label{gw_model}

\begin{figure}[bht]
\includegraphics[width=4.83in,height=4.38in]{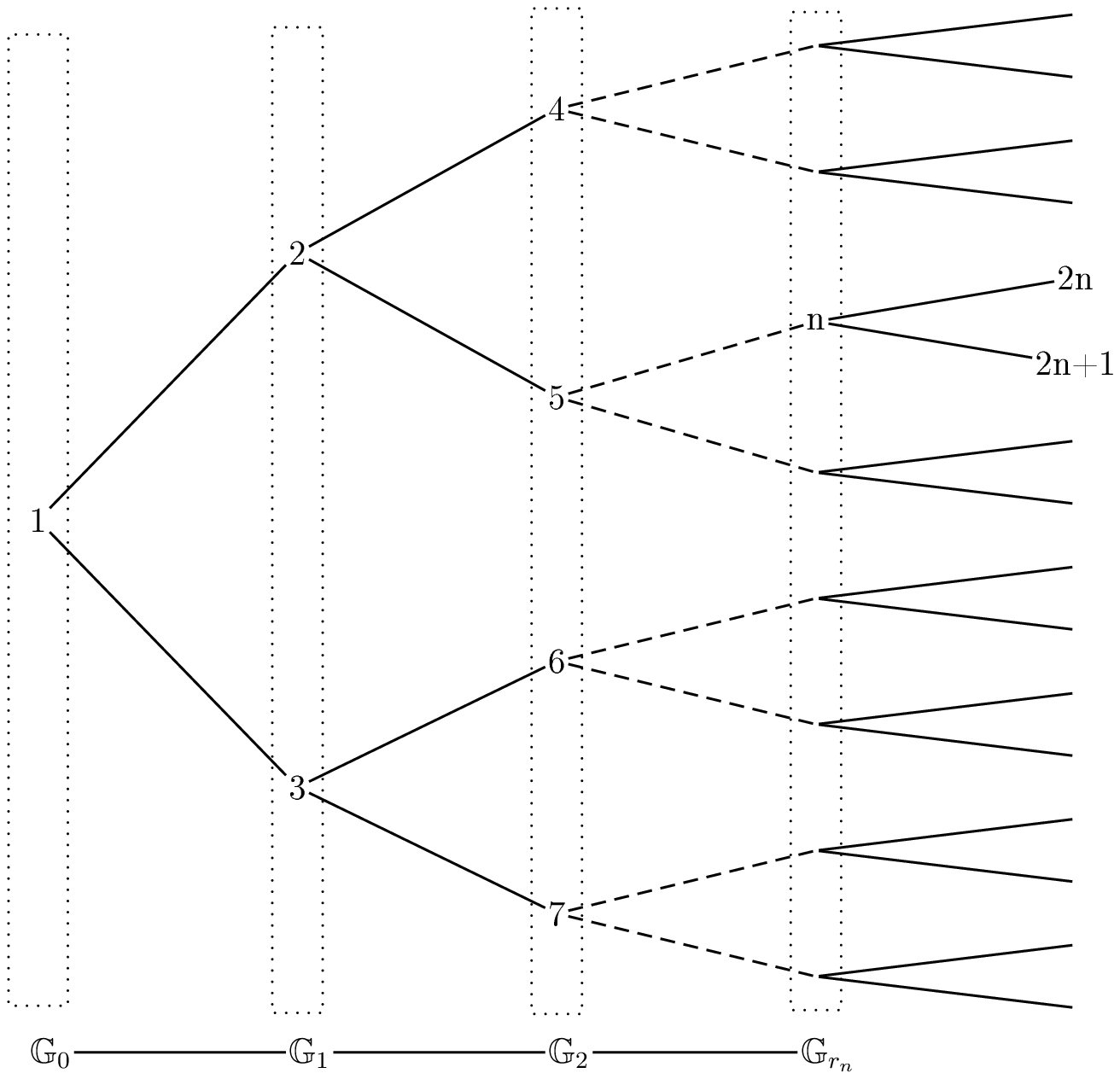}
\caption{The binary tree $\mathbb{T}$}\label{Fi:gw_arbrebinaire}
\end{figure}

Let $\mathbb{T}$ be a binary regular tree in which each vertex is
seen as a positive integer different from 0, see Figure
\ref{Fi:gw_arbrebinaire}. For $r\in\mathbb{N}$, let
\[
\mathbb{G}_{r}=\Big\{2^{r},2^{r}+1,\cdots,2^{r+1}-1\Big\}, \quad
\mathbb{T}_{r}=\bigcup\limits_{q=0}^{r}\mathbb{G}_{q},
\]
which denote respectively the $r$-th column and the first $(r+1)$
columns of the tree.
Then, the cardinality $|\mathbb{G}_{r}|$ of $\mathbb{G}_{r}$ is
$2^{r}$ and that of $\mathbb{T}_{r}$ is
$|\mathbb{T}_{r}|=2^{r+1}-1$. A column of a given integer $n$ is
$\mathbb{G}_{r_{n}}$ with $r_{n}=\lfloor\log_{2}n\rfloor$, where
$\lfloor x\rfloor$ denotes the integer part of the real number $x$.

The genealogy of the cells is described by this tree. In the sequel
we will thus see $\mathbb{T}$ as a given population. Then the vertex
$n$, the column $\mathbb{G}_{r}$ and the first $(r+1)$ columns
$\mathbb{T}_{r}$ designate respectively individual $n$, the $r$-th
generation and the first $(r+1)$ generations. The initial individual
is denoted $1$. The model proposed by Delmas and Marsalle
\cite{DelMar} is defined as follows. The growth rate of cell $n$ is
$X_{n}.$
\begin{itemize}
\item With probability $p_{1,0},$ $n$ gives birth to two cells $2n$
and $2n+1$ with both divide. The growth rate of the daughters
$X_{2n}$ and $X_{2n+1}$ are then linked to the mother's one through
auto-regressive equations (\ref{gw_bar11}).
\item With probability $p_{0},$ only the new pole $2n$ divides. Its
growth rate $X_{2n}$ is linked to its mother's one $X_{n}$ through
the relation
\begin{equation}\label{eq_new}
X_{2n}=\alpha_{0}'X_{n}+\beta_{0}'+\vep_{2n}',
\end{equation}
where $\alpha_{0}'\in(-1,1),$ $\beta_{0}'\in\R$ and $(\vep_{2n}',
n\in\T)$ is a sequence of independent centered Gaussian random
variables with variance $\sigma_{0}^{2}>0.$
\item With probability $p_{1},$ only the old pole $2n+1$ divides. Its
growth rate $X_{2n+1}$ is linked to its mother's one $X_{n}$ through
the relation
\begin{equation}\label{eq_old}
X_{2n+1}=\alpha_{1}'X_{n}+\beta_{1}'+\vep_{2n+1}',
\end{equation}
where $\alpha_{1}'\in(-1,1),$ $\beta_{1}'\in\R$ and $(\vep_{2n+1}',
n\in\T)$ is a sequence of independent centered Gaussian random
variables with variance $\sigma_{1}^{2}>0.$
\item With probability $1-p_{1,0}-p_{1}-p_{0},$ which is
non-negative, $n$ gives birth to two cells which do not divide.
\item The sequences $((\vep_{2n},\vep_{2n+1}),n\in\T),$
$(\vep_{2n}',n\in\T)$ and $(\vep_{2n+1}',n\in\T)$ are independent.
\end{itemize}
The process $(X_{n})$ described above is a typical example of BMC on
GW tree. In \cite{DesaporGegMar}, this process is called bifurcating
autoregressive process (BAR) with missing data. It is an extension
of bifurcating autoregressive process when the data have a binary GW
tree structure, see figure 2 for example of binary GW tree. Indeed,
one can assume that the cells which do not divide and those which do
not exist are missing or dead.

In \cite{DelMar}, Delmas and Marsalle using their results for BMC on
GW tree, gave laws of large numbers and central limit theorem for
the maximum likelihood estimator of the parameter
\begin{equation}\label{theta}
\theta =
(\alpha_{0},\beta_{0},\alpha_{1},\beta_{1},\alpha_{0}',\beta_{0}',\alpha_{1}',\beta_{1}').
\end{equation}
In this paper, we will give deviation inequalities for the least
squares estimator of the parameter $\theta,$ in case the noise
sequence and the initial state $X_{1}$ take their values in a
compact set. Note that this implies that the BAR process with
missing data describes above also take their values in compact set.
These deviation inequalities are important for a rigorous non
asymptotic statistical study. Indeed, when the sample size is
insufficient to apply limit theorems, they allow for example to
estimate the errors in the estimation of unknown parameters.
Furthermore, these inequalities allow to get a rate of convergence
in the laws of large numbers, and this permit, for example, to build
non-asymptotic confidence intervals.


We are now going to give a rigorous definition of BMC on GW tree. We
refer to \cite{DelMar} for more details.

\subsection{Definitions}

For an individual $n\in \mathbb{T}$, we are interested in the
quantity $X_{n}$ (it may be the weight, the growth rate,$\cdots$)
with values in the metric space $S$ endowed with its Borel
$\sigma$-field $\mathcal{S}$.
\begin{definition}[$\mathbb{T}$-transition probability, see (\cite{Guyon})]
We call $\mathbb{T}$-transition probability any mappings $P: S\times
\mathcal{S}^{2}\rightarrow [0,1]$ such that
\begin{itemize}
\item $P(.,A)$ is measurable for all $A\in \mathcal{S}^{2}$,

\item $P(x,.)$ is a probability measure on $(S^{2},\mathcal{S}^{2})$
for all $x\in S$.
\end{itemize}
\end{definition}

For $p\geq 1$, we denote by $\mathcal{B}(S^{p})$(resp.
$\mathcal{B}_{b}(S^{p})$, $\mathcal{C}(S^{p})$,
$\mathcal{C}_{b}(S^{p})$) the set of all
$\mathcal{S}^{p}$-measurable (resp. $\mathcal{S}^{p}$-measurable and
bounded, continuous, continuous and bounded) mapping $f:
S^{p}\rightarrow \mathbb{R}$. For $f\in \mathcal{B}(S^{3})$, when it
is defined, we denote by $Pf \in \mathcal{B}(S)$ the function
\[x\mapsto Pf(x)=\int_{S^{2}}f(x,y,z)P(x,dy,dz).\]
\begin{definition}[Bifurcating Markov Chains, see (\cite{Guyon})]
Let $(X_{n}, n\in \mathbb{T})$ be a family of $S$-valued random
variables defined on a filtered probability space $(\Omega,
\mathcal{F}, (\mathcal{F}_{r}, r\in\mathbb{N}), \mathbb{P})$. Let
$\nu$ be a probability on $(S, \mathcal{S})$ and $P$ be a
$\mathbb{T}$-transition probability. We say that $(X_{n}, n\in
\mathbb{T})$ is a $(\mathcal{F}_{r})$-bifurcating Markov chain with
initial distribution $\nu$ and $\mathbb{T}$-transition probability
$P$ if
\begin{itemize}
\item $X_{n}$ is $\mathcal{F}_{r_{n}}$-measurable for all $n\in
\mathbb{T}$,
\item $\mathcal{L}(X_{1})=\nu$,
\item for all $r\in \mathbb{N}$ and for all family $(f_{n}, n\in
\mathbb{G}_{r})\subseteq \mathcal{B}_{b}(S^{3})$
\[\mathbb{E}\left[\prod\limits_{n\in
\mathbb{G}_{r}}f_{n}(X_{n},X_{2n},X_{2n+1})\Big|
\mathcal{F}_{r}\right]=\prod\limits_{n\in
\mathbb{G}_{r}}Pf_{n}(X_{n}).\]
\end{itemize}
\end{definition}

Now, we add a cemetery point to $S$, $\partial.$ Let
$\bar{S}=S\cup\{\partial\},$ and $\bar{\mathcal{S}}$ be the
$\sigma-$field generated by $\mathcal{S}$ and $\{\partial\}.$ In the
previous biological framework, $S$ corresponds to the state space of
the quantities related to living cells, and $\partial$ is the
default value for dead cells. Let $P^*$ be a $\T$-transition
probability defined on $\bar{S} \times \bar{{\mathcal{S}}}$ such
that
\begin{equation}\label{P_ast}
P^{*}(\partial,\{(\partial,\partial)\})=1.
\end{equation}
In the previous biological framework, (\ref{P_ast}) means that no
dead cell can give birth to a living cell. We denote by $P_{0}^{*}$
and $P_{1}^{*}$ the restriction of the first and the second marginal
of $P^*$ to $S$, that is:
\begin{equation*}
P_{0}^{*}=P^{*}\left(\cdot,\left(\cdot\bigcap S\right)\times
\bar{S}\right) \quad \text{and} \quad
P_{1}^{*}=P^{*}\left(\cdot,\bar{S}\times\left(\cdot\bigcap
S\right)\right).
\end{equation*}

\begin{definition}[BMC on GW tree, see \cite{DelMar}]
Let $X=(X_{n},n\in\T)$ be a $P^{*}$-BMC on $(\bar{S},\bar{\SS}),$
with $P^{*}$ satisfying (\ref{P_ast}). We call $(X_{n},n\in\T^{*}),$
with $\T^{*} = \{n\in\T:X_{n}\neq \partial\},$ a BMC on GW tree. The
$P^{*}$-BMC is said spatially homogeneous if
$p_{1,0}=P^{*}(x,S\times S),$ $p_{0}=P^{*}(x,S\times \{\partial\}),$
and $p_{1}=P^{*}(x,\{\partial\}\times S)$ do not depend on $x\in S.$
A spatially homogeneous $P^{*}$-BMC is said super-critical if $m>1,$
where $m=2p_{1,0}+p_{1}+p_{0}.$
\end{definition}

We denote by $(Y_{n},n\in\N)$ the Markov chain on $S$ with
$Y_{0}=X_{1}$ and transition probability
$Q=\frac{1}{m}(P_{0}^{*}+P_{1}^{*}).$

\begin{remark}

\begin{itemize}
\item The name BMC on GW tree comes from the fact that condition
(\ref{P_ast}) and spatial homogeneity imply that $\T^*$ is a GW
tree.
\item All through this work, we shall assume that the $P^*$-BMC is
super-critical.
\end{itemize}
\end{remark}

\begin{figure}[bht]
\includegraphics[width=3.83in,height=3.38in]{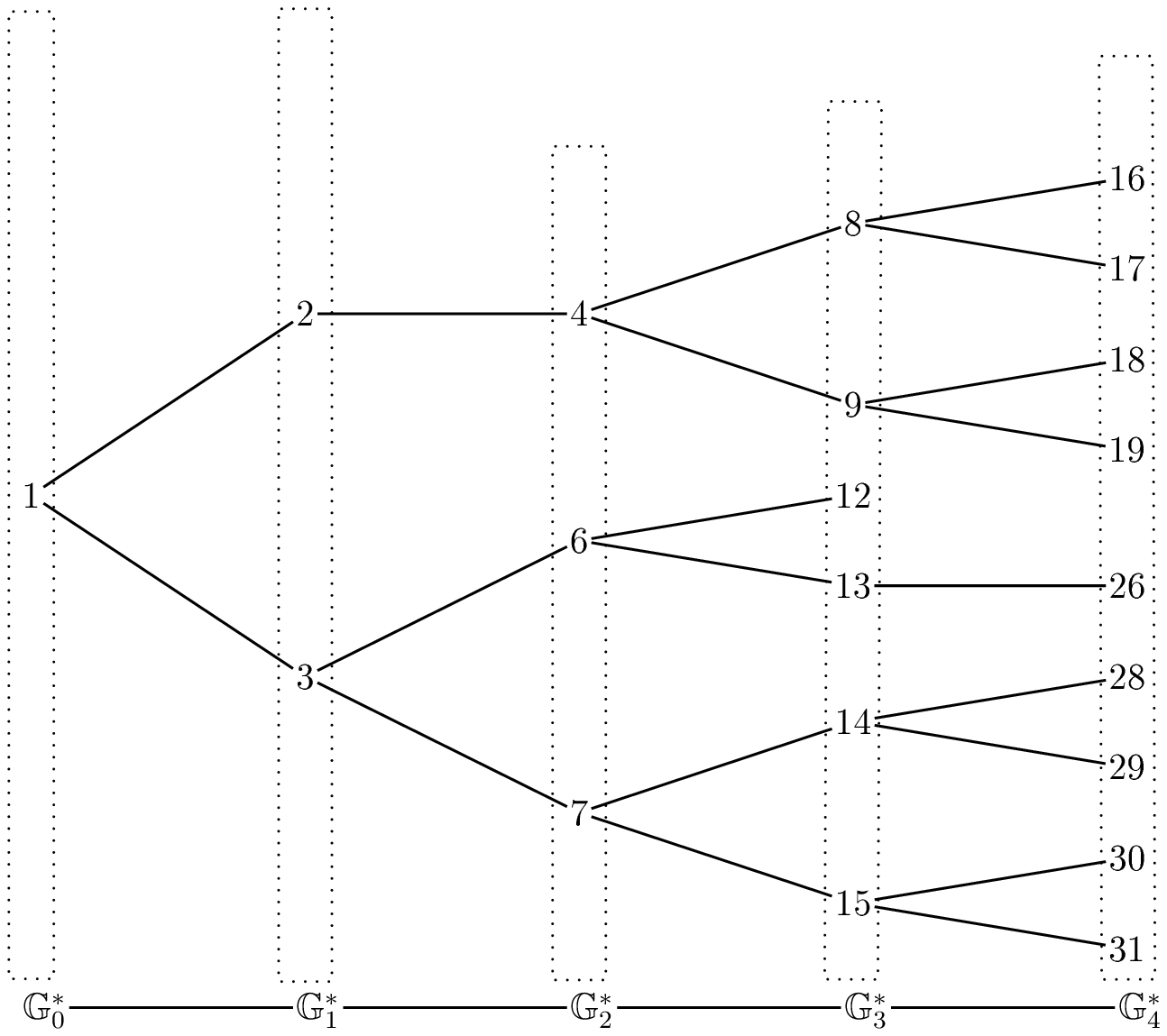}
\caption{A binary GW tree up to the 4 th generation. \emph{In this
tree, individual 1 gives birth to two individuals which both divide,
this happen with probability $p_{1,0}.$ Individual 2 gives birth to
two individuals which only one (the new pole) divides, this happen
which probability $p_{0}.$ Individual 12 gives birth to two
individuals which do not divide, this happen with probability
$1-p_{1,0}-p_{0}-p_{1}.$ }} \label{Fi:gwtree}
\end{figure}

Now, for any subset $J\subset\T,$ let
\begin{equation*}
J^{*}=J\cap\T^{*}=\{j\in J: X_{j}\neq\partial\}
\end{equation*}
be the subset of living cells among $J,$ and $|J|$ be the cardinal
of $J.$ The process $\left(|\G_{k}^{*}|,k\in\N\right),$ is a GW
process with the reproduction generating function
\begin{equation*}
\psi(z)=(1-p_{0}-p_{1}-p_{1,0}) + (p_{0} + p_{1})z + p_{1,0}z^{2},
\end{equation*}
and the average number of daughters alive is $m.$ It is known, see
e.g \cite{AtNe72}, that $m^{-k}|\G_{k}^{*}|$ converges in
probability to a non-negative random variable $W$. Moreover,
$\P(W>0)=1$ iff there is no extinction. We have for all $r\geq 0,$
\begin{equation}\label{t_r}
\E\left[|\G_{r}^{*}|\right]=m^{r} \quad \text{and} \quad
\E\left[|\T_{r}^{*}|\right] =
\sum\limits_{q=0}^{r}\E\left[|\G_{q}^{*}|\right] =
\frac{m^{r+1}-1}{m-1} := t_{r}.
\end{equation}
It is known, see \cite{DelMar}, that $t_{r}^{-1}|\T_{r}^{*}|$
converges in probability to $W$ as well.

For $i\in\T,$ set $\Delta_{i}=(X_{i},X_{2i},X_{2i+1})$ the
mother-daughters quantities of interest. For a finite subset
$J\subset\T,$ we set
\begin{equation}\label{sample_sum}
M_{J}(f)=\begin{cases} \sum\limits_{i\in J}f(X_{i}) \quad \text{for
$f\in\BB(\bar{S}),$} \\ \sum\limits_{i\in J}f(\Delta_{i}) \quad
\text{for $f\in \BB(\bar{S}^{3}),$}\end{cases}
\end{equation}
with the convention that a sum over an empty set is null. We also
define the following two averages of $f$ over $J$
\begin{equation}\label{moy}
\overline{M}_{J}(f) = \frac{1}{|J|}M_{J}(f) \quad \text{if $|J|>0$}
\quad \text{and} \quad
\widetilde{M}_{J}(f)=\frac{1}{\E\left[|J|\right]}M_{J}(f) \quad
\text{if $\E\left[|J|\right]>0.$}
\end{equation}
Limit theorems for averages (\ref{moy}) have been studied in
\cite{DelMar} for $J = \G_{n}^{*}$ and $J=\T_{n}^{*},$ as $n$ goes
to infinity. Under uniform geometric ergodicity assumption for $Q$,
we will establish in this paper deviation inequalities for those
averages. These deviation inequalities will allow to highlight three
regimes for the speed of convergence of above averages, thus showing
a competition between the ergodicity of the embedded Markov chain
$\displaystyle (Y_{n}, n\in\mathbb{N})$ and the size of the binary
Galton-Watson tree. This new phenomenon is not observed in the
asymptotic study of Delmas and Marsalle \cite{DelMar}. Notice that
deviation inequalities were already studied in the no death case
\cite{BDG11}, that is $m=2.$ We will follow essentially the same
approach that the latter paper for the proofs of our results.
However, we will introduce some modifications on those proofs in
order to take into account the randomness of index set, and we will
make use of the theory of large deviation for branching processes
\cite{At94}. Let us also mention \cite{BiDj12}, where the authors
establish deviation inequalities for estimators of parameters of the
$p$-order bifurcating autoregressive process.

The rest of paper is organized as follows. In section \ref{gw_main},
we states our main results, that is deviation inequalities for
averages (\ref{moy}), for $J = \G_{n}^{*}$ and $J=\T_{n}^{*}$. This
will be done under uniform geometric ergodicity assumption for $Q$,
and suitable assumptions on the binary GW tree. In section
\ref{gw_appli}, we will focus in particular on the first order
bifurcating autoregressive process with missing data described in
section \ref{gw_model}. Section \ref{gw_proofs} is dedicated to the
proofs of our results.

\section{Main results}\label{gw_main}
We consider the following hypothesis:
\begin{enumerate}
\item [(\bf H1):] There exists a probability measure $\mu$ on $(S,\SS)$ such that for all
$f\in\BB_{b}(S)$ with $\langle\mu,f\rangle=0,$ there is $c>0$ such
that for all $k\in\N$ and for all $x\in S,$ $|Q^{k}f(x)|\leq
c\alpha^{k}.$
\item [{\bf (H2):}] $m>\sqrt{2}$.
\item [{\bf (H3):}] $p_{1,0} + p_{0} + p_{1} = 1$, where $p_{1,0}$, $p_{0}$ and
$p_{1}$ are defined in section \ref{gw_model}.
\end{enumerate}

\begin{remark}
Hypothesis {\bf (H1)} implies that the Markov chain $Y$ is ergodic,
that is for all $f\in\CC_{b}(S)$ and for all $x\in S,$
$\lim\limits_{k\rightarrow\infty}\E_{x}[f(Y_{k})]=\langle\mu,f\rangle.$
Assuming hypothesis \textbf{(H3)} means that we work conditionally
to the non-extinction. Note that this is consistent with the study
of E. Coli.

\medskip

Hypothesis ({\textbf{H2}}) comes from our calculations. indeed, in
order to get relevant inequalities, i.e. inequalities for which the
upper bound goes to zero as the sample size increases, we have to
assume that $m>\sqrt{2}$. However, our deviation inequalities also
work for $m\leq\sqrt{2}$, but they are not relevant for this case.
To get relevant deviation inequalities for $m\leq\sqrt{2}$ is still
an open problem that we will pursue in an other work.
\end{remark}
In the sequel, $\H_{r}$ will denote one of the set $\G_{r}$ or
$\T_{r}.$ We set $h_{r}=(m^{2}/2)^{r}$ if $\H_{r}=\G_{r}$ and
$h_{r}=(m^{2}/2)^{r+1}$ if $\H_{r}=\T_{r}$. We can now state our
main results. Notice that any function $f$ defined on $S$ is
extended to $\bar{S}$ by setting $f(\partial)=0.$

\begin{thm}\label{thm:gw_main1}
Under hypothesis {\bf (H1)} and {\bf (H2)}, let $f\in \BB_{b}(S)$
such that $\langle\mu,f\rangle=0.$ Then we have for all $\delta>0$:

\begin{itemize}
\item if $m\alpha<1$, then $\forall r\in\N$,
\begin{equation*}
\mathbb{P}\Big(\widetilde{M}_{\H_{r}^{*}}(f)>\delta\Big)\leq
\exp\left(c''\delta\right)\exp\left(-c'\delta^{2}h_{r}\right);
\end{equation*}
\item if $m\alpha= 1$, then for $\H_{r}=\G_{r}$ and $\forall r\in \N$,
\begin{equation*}
\mathbb{P}\Big(\widetilde{M}_{\H_{r}^{*}}(f)>\delta\Big)\leq
\exp\left(c''\delta\right)\exp\left(-c'\delta^{2}(m^{2}/2)^{r}\right);
\end{equation*}
\item if $m\alpha = 1$, then for $\H_{r}=\T_{r}$ and $\forall r\in\N$,
\begin{equation*}
\mathbb{P}\Big(\widetilde{M}_{\H_{r}^{*}}(f)>\delta\Big)\leq
\exp\left(c''\delta(r+1)\right)\exp\left(-c'\delta^{2}(m^{2}/2)^{r+1}\right);
\end{equation*}
\item if $1<m\alpha<\sqrt{2}$, then $\forall r \in\N$ such that $r>r_0$,
\begin{equation*}
\mathbb{P}\Big(\widetilde{M}_{\H_{r}^{*}}(f)>\delta\Big)\leq
\exp\left(-c'\delta^{2}h_{r}\right);
\end{equation*}
\item if $m\alpha=\sqrt{2}$, then $\forall r\in \N$  such that $r>r_0$,
\begin{equation*}
\mathbb{P}\Big(\widetilde{M}_{\H_{r}^{*}}(f)>\delta\Big)\leq
\exp\left(-\frac{c'\delta^{2}h_{r}}{r}\right);
\end{equation*}
\item if $m\alpha>\sqrt{2}$, then $\forall r\in \N^{*}$ such that $r>r_0$,
\begin{equation*}
\mathbb{P}\Big(\widetilde{M}_{\H_{r}^{*}}(f)>\delta\Big)\leq
\exp\left(-\frac{c'\delta^{2}}{\alpha^{2r}}\right);
\end{equation*}
\end{itemize}
where,
\begin{itemize}
\item $r_0:= \log\left(\frac{\delta}{c_{0}}\right)/\log(\alpha) -
k_{0},$ with $k_{0}\in\{0,1\}$,
\item $c_{0}$, $c'$ and $c''$ are positive constants which depend on $\alpha$, $m$, and $c$ and
may differ line by line.
\end{itemize}
\end{thm}

\begin{thm}\label{thm:gw_main11}
Under hypothesis {\bf (H1)}-{\bf (H3)}, we have for all $f\in
\BB_{b}(S)$ such that $\<\mu,f\>\neq0$ and for all $\delta>0$:
\begin{itemize}
\item if $m\alpha<1$, then $\forall r\in\N$,
\begin{equation*}
\mathbb{P}\Big(\widetilde{M}_{\H_{r}^{*}}(f)-\<\mu,f\>W>
\delta\Big)\leq
\exp\left(c''\delta\right)\exp\left(-c'\delta^{2}h_{r}\right) +
A_{r};
\end{equation*}
\item if $m\alpha= 1$, then for $\H_{r}=\G_{r}$ and $\forall r\in \N$,
\begin{equation*}
\mathbb{P}\Big(\widetilde{M}_{\H_{r}^{*}}(f)-\<\mu,f\>W>\delta\Big)\leq
\exp\left(c''\delta\right)\exp\left(-c'\delta^{2}(m^{2}/2)^{r}\right)
+ A_{r};
\end{equation*}
\item if $m\alpha = 1$, then for $\H_{r}=\T_{r}$ and $\forall r\in\N$,
\begin{equation*}
\mathbb{P}\Big(\widetilde{M}_{\H_{r}^{*}}(f)-\<\mu,f\>W>\delta\Big)\leq
\exp\left(c''\delta(r+1)\right)\exp\left(-c'\delta^{2}(m^{2}/2)^{r+1}\right)
+ A_{r};
\end{equation*}
\item if $1<m\alpha<\sqrt{2}$, then $\forall r \in\N$ such that $r>r_0$,
\begin{equation*}
\mathbb{P}\Big(\widetilde{M}_{\H_{r}^{*}}(f)-\<\mu,f\>W>\delta\Big)\leq
\exp\left(-c'\delta^{2}h_{r}\right) + A_{r};
\end{equation*}
\item if $m\alpha=\sqrt{2}$, then $\forall r\in \N$  such that $r>r_0$,
\begin{equation*}
\mathbb{P}\Big(\widetilde{M}_{\H_{r}^{*}}(f)-\<\mu,f\>W>\delta\Big)\leq
\exp\left(-\frac{c'\delta^{2}h_{r}}{r}\right) + A_{r};
\end{equation*}
\item if $m\alpha>\sqrt{2}$, then $\forall r\in \N^{*}$ such that $r>r_0$,
\begin{equation*}
\mathbb{P}\Big(\widetilde{M}_{\H_{r}^{*}}(f)-\<\mu,f\>W>\delta\Big)\leq
\exp\left(-\frac{c'\delta^{2}}{\alpha^{2r}}\right) + A_{r};
\end{equation*}
\end{itemize}
where,
\begin{itemize}
\item for all $r\in\N,$ \begin{equation*}
A_{r} = \begin{cases}
c'\exp\left(-c''\delta^{2/3}(m^{1/3})^{r}\right) \hspace{3.5cm}
\text{if
$\H_{r}=\G_{r}$} \\
\exp\left(c'\delta^{2/3}\right)\exp\left(-c''\delta^{2/3}\left(t_{r}/(r+1)^{2}\right)^{1/3}\right)
\quad \text{if $\H_{r}=\T_{r},$}
\end{cases}
\end{equation*}
\item $r_0:= \log\left(\frac{\delta}{c_{0}}\right)/\log(\alpha) -
k_{0},$ with $k_{0}\in\{0,1\},$
\item $c_{0},$ $c'$ and $c''$
are positive constants which depend on $\alpha$, $m$, and $c$ and
may differ line by line.
\end{itemize}
\end{thm}

\begin{remark}
For $\<\mu,f\>=0$ (Theorem \ref{thm:gw_main1}), there is no
additional term in the deviation of average
$\widetilde{M}_{\H_{r}^{*}}(f)$. While in Theorem
\ref{thm:gw_main11} there is an additional term $A_{r}$ which
appears. This term is a contribution of the binary Galton-Watson
tree on the deviation of average $\widetilde{M}_{\H_{r}^{*}}(f)$
with respect to $\<\mu,f\>W$. This explain why we need additional
hypothesis (\textbf{H3}) in Theorem \ref{thm:gw_main11}, because we
have to deal with the deviation inequalities for Galton-Watson
processes.
\end{remark}
The next results can be seen as a consequence of the previous
results.
\begin{thm}\label{thm:gw_main2}
We assume that hypothesis {\bf (H1)-(H3)} are satisfied. Let
$f\in\BB_{b}(S).$ For all $\delta>0,$ for all $a>0$ and for all
$b>0$ such that $b<a/(\delta+1),$ we have
\begin{itemize}
\item if $m\alpha<1$, then $\forall r\in\N$,
\begin{equation*}
\mathbb{P}\Big(\overline{M}_{\H_{r}^{*}}(f)-\<\mu,f\>
>\delta\Big|W\geq a\Big)\leq
\exp\left(c''\delta b\right)\exp\left(-c'(\delta b)^{2}h_{r}\right)
+ A_{r};
\end{equation*}
\item if $m\alpha= 1$, then for $\H_{r}=\G_{r}$ and $\forall r\in \N$,
\begin{equation*}
\mathbb{P}\Big(\overline{M}_{\H_{r}^{*}}(f)-\<\mu,f\>
>\delta\Big|W\geq a\Big)\leq
\exp\left(c''\delta b\right)\exp\left(-c'(\delta
b)^{2}(m^{2}/2)^{r}\right) + A_{r};
\end{equation*}
\item if $m\alpha = 1$, then for $\H_{r}=\T_{r}$ and $\forall r\in\N$,
\begin{equation*}
\mathbb{P}\Big(\overline{M}_{\H_{r}^{*}}(f)-\<\mu,f\>
>\delta\Big|W\geq a\Big)\leq
\exp\left(c''\delta b(r+1)\right)\exp\left(-c'(\delta
b)^{2}(m^{2}/2)^{r+1}\right) + A_{r};
\end{equation*}
\item if $1<m\alpha<\sqrt{2}$, then $\forall r \in\N$ such that $r>r_0$,
\begin{equation*}
\mathbb{P}\Big(\overline{M}_{\H_{r}^{*}}(f)-\<\mu,f\>
>\delta\Big|W\geq a\Big)\leq
\exp\left(-c'(\delta b)^{2}h_{r}\right) + A_{r};
\end{equation*}
\item if $m\alpha=\sqrt{2}$, then $\forall r\in \N$  such that $r>r_0$,
\begin{equation*}
\mathbb{P}\Big(\overline{M}_{\H_{r}^{*}}(f)-\<\mu,f\>
>\delta\Big|W\geq a\Big)\leq
\exp\left(-\frac{c'(\delta b)^{2}h_{r}}{r}\right) + A_{r};
\end{equation*}
\item if $m\alpha>\sqrt{2}$, then $\forall r\in \N^{*}$ such that $r>r_0$,
\begin{equation*}
\mathbb{P}\Big(\overline{M}_{\H_{r}^{*}}(f)-\<\mu,f\>
>\delta\Big|W\geq a\Big)\leq
\exp\left(-\frac{c'(\delta b)^{2}}{\alpha^{2r}}\right) + A_{r};
\end{equation*}
\end{itemize}
where,
\begin{itemize}
\item for all $r\in\N,$ \begin{equation*}
A_{r} = \begin{cases} c'\exp\left(-c''(\delta
b)^{2/3}(m^{1/3})^{r}\right) \hspace{4cm} \text{if
$\H_{r}=\G_{r}$} \\
\exp\left(c'(\delta b)^{2/3}\right)\exp\left(-c''(\delta
b)^{2/3}\left(t_{r}/(r+1)^{2}\right)^{1/3}\right) \quad \text{if
$\H_{r}=\T_{r},$}
\end{cases}
\end{equation*}
\item $r_0:= \log\left(\frac{\delta b}{c_{0}}\right)/\log(\alpha) -
k_{0},$ with $k_{0}\in\{0,1\},$
\item $c_{0},$ $c'$ and $c''$
are positive constants which depend on $\alpha$, $m,$ $a,$ and $c$,
and may differ line by line.
\end{itemize}
\end{thm}
We have the following extension of above theorems when $f$ does not
only depend on an individual $X_{i},$ but on the mother-daughters
triangle $\Delta_{i}.$
\begin{thm}\label{thm:gw_main3}
Let $f\in\BB_{b}(S^{3}).$ If $\<\mu,P^{*}f\>=0,$ then, under
hypothesis \textbf{(H1)} and \textbf{(H2)}, we have deviation
inequalities of Theorem \ref{thm:gw_main1} for
$\widetilde{M}_{\H_{r}^{*}}(f).$ If $\<\mu,P^{*}f\>\neq0,$ under
additional hypothesis \textbf{(H3)}, we have deviation inequalities
of Theorem \ref{thm:gw_main11} for $\widetilde{M}_{\H_{r}^{*}}(f) -
\<\mu,P^{*}f\>W$ and of Theorem \ref{thm:gw_main2} for
$\overline{M}_{\H_{r}^{*}}(f) - \<\mu,P^{*}f\>$.
\end{thm}

\begin{remark}
Let us stress that by tedious, but straightforward calculations, the
constants which appear in the previous inequalities can be made
explicit.
\end{remark}

Let us recall the following definition.

\begin{definition}\label{exponentialdefini}
Let $(E,d)$ be a metric space. Let $(Z_{n})$ be a sequence of random
variables valued in $E$, $Z$ be a random variable valued in $E$ and
$(v_{n})$ be a rate. We  say that $Z_n$ converges
$v_n$-superexponentially  fast in probability to $Z$ if for all
$\delta>0$,
\begin{equation*}
\limsup_{n\rightarrow
\infty}\frac{1}{v_n}\log\P(d(Z_n,Z)>\delta)=-\infty.
\end{equation*}
This ``exponential convergence'' with speed $v_n$ will be shortened
as
\begin{equation*}
Z_n \overset {\rm superexp}{\underset{v_n}{\longrightarrow}} Z.
\end{equation*}
\end{definition}

\begin{remark}
Let $(b_{n})$ be a sequence of increasing positive real numbers such
that
\begin{equation*}
b_{n} \rightarrow +\infty
\end{equation*}
and
\begin{itemize}
\item if $m\alpha<\sqrt{2}$, the sequence $(b_{n})$ is such that $\displaystyle b_{n}/\sqrt{n}\longrightarrow 0$,

\item if $m\alpha=\sqrt{2}$, the sequence $(b_{n})$ is such that $\displaystyle (b_{n}\sqrt{\log n})/\sqrt{n}\longrightarrow 0$,

\item  if $m\alpha>\sqrt{2}$, the sequence $(b_{n})$ is such that
$\displaystyle b_{n}\alpha^{\log n/\log (m^{2}/2)}\longrightarrow
0$.
\end{itemize}

From the previous deviations inequalities, we can deduce easily that
\begin{equation*}
\widetilde{M}_{\mathbb{H}_{r}^{*}}(f) \overset {\rm
superexp}{\underset{b_{\lfloor h_{r}\rfloor}^{2}}{\longrightarrow}}
0 \quad \text{if $\langle \mu,f\rangle = 0$},
\end{equation*}
and if $\langle\mu,f\rangle\neq0$, we have for $m<2^{3/5}$
\begin{equation*}
\widetilde{M}_{\mathbb{H}_{r}^{*}}(f) \overset {\rm
superexp}{\underset{b_{\lfloor h_{r}\rfloor}^{2}}{\longrightarrow}}
\langle\mu,f\rangle W,
\end{equation*}
and $\forall a>0$,
\begin{equation*}
\limsup_{r\rightarrow+\infty} \frac{1}{b_{\lfloor h_{r}\rfloor}^{2}}
\log\P\left(\left|\overline{M}_{\mathbb{H}_{r}^{*}}(f) -
\langle\mu,f\rangle\right|>\delta\big|W\geq a\right) = -\infty.
\end{equation*}
So, for the exponential convergence of averages
$\widetilde{M}_{\mathbb{H}_{r}^{*}}(f)$ and
$\overline{M}_{\mathbb{H}_{r}^{*}}(f)$, there are three regimes
according to the value of $m\alpha$ compared to $\sqrt{2}$. This
phenomenon is not observed in the limit theorems of Delmas and
Marsalle \cite{DelMar}. However, a similar phenomenon was observed
recently by Adamczak and Mi{\l}o\'s for the central limit theorem of
branching particle system \cite{AM}.

\medskip

So, our deviations inequalities highlight a competition between the
ergodicity of the embedded Markov chain with transition probability
$Q$ and the Galton-Watson binary tree.
\end{remark}

\section{Application:First order bifurcating autoregressive
processes with missing data}\label{gw_appli}

We consider the asymmetric auto-regressive processes given in
section \ref{gw_model}. Notice that the process $(X_{i},i\in\T)$
defined in section \ref{gw_model}, with the convention that
$X_{i}=\partial$ if the cell $i$ is missing, is a spatially
homogeneous BMC on a GW tree. We will assume that
$2p_{1,0}+p_{1}+p_{0}>\sqrt{2}$. This implies in particular that the
BMC on GW is super-critical. We will also assume that the noise
sequences $((\vep_{2n},\vep_{2n+1}),n\in\T),$ $(\vep_{2n}',n\in\T)$
and $(\vep_{2n+1}',n\in\T)$, and the initial state $X_{1}$ take
their values in a compact set. The latter implies that the process
$(X_{i},i\in\T)$ is bounded. We denote by $S$ the state space of
$(X_{i},i\in\T)$. We assume without loss of generality that $S$ is a
compact subset of $\R$.

Let $\T_{n}^{0,1}$ be the subset of cells in $\T_{n}^{*}$ with two
living daughters, $\T_{n}^{0}$ (resp.$\T_{n}^{1}$) be the set of
cells of $\T_{n}^{*}$ with only the new (resp. old) pole daughter
alive:

\hspace*{2cm}$ \T_{n}^{1,0} = \left\{i\in\T_{n}^{*}: \Delta_{i} \in
S^{3}\right\},  \, \, \, \, \, \, \, \, \, \, \T_{n}^{0} =
\{i\in\T_{n}^{*}: \Delta_{i}\in S^{2}\times \{\partial\}\} $  \, \,
\, and \, \, \, \, \, \, \, \, \, \, \, \, \, \, \hspace*{2cm}
$\T_{n}^{1} = \{i\in\T_{n}^{*}:\Delta_{i}\in
S\times\{\partial\}\times S\}.$

We compute the least-squares estimator (LSE)
\begin{equation*}
\widehat{\theta}_{n} = (\widehat{\alpha}_{0}^{n},
\widehat{\beta}_{0}^{n}, \widehat{\alpha}_{1}^{n},
\widehat{\beta}_{1}^{n}, \widehat{\alpha}_{0}^{'n},
\widehat{\beta}_{0}^{'n}, \widehat{\alpha}_{1}^{'n},
\widehat{\beta}_{1}^{'n})
\end{equation*}
of $\theta$ given by (\ref{theta}), based on the observation of a
sub-tree $\T_{n+1}^{*}.$ Consequently, we obviously have for
$\eta\in\{0,1\},$

\vspace{0.25cm}

$\widehat{\alpha}_{\eta}^{n}  =
\frac{\displaystyle{|\T_{n}^{1,0}|^{-1}
     \sum_{i \in \T_{n}^{1,0}}X_{i}X_{2i+\eta} - \left(|\T_{n}^{1,0}|^{-1} \sum_{i
       \in \T_{n}^{1,0}}X_i\right)\left(|\T_{n}^{1,0}|^{-1}\sum_{i \in \T_{n}^{1,0}}
     X_{2i+\eta}\right)}}{\displaystyle{|\T_{n}^{1,0}|^{-1} \sum_{i \in
       \T_{n}^{1,0}}X_i^2 - \left(|\T_{n}^{1,0}|^{-1}\sum_{i \in
       \T_{n}^{1,0}}X_i\right)^2}}, \vspace{0.022cm} \\ \widehat{\beta}_{\eta}^{n}  =  \displaystyle{|\T_n^{1,0}|^{-1} \sum_{i \in
     \T_{n}^{1,0}} X_{2i+\eta} - \widehat{\alpha}_{\eta}^{n} |\T_{n}^{1,0}|^{-1}
   \sum_{i \in \T_n^{1,0}} X_i},\\ \widehat{\alpha}_{\eta}'^{n}  =  \frac{\displaystyle{|\T_{n}^{\eta}|^{-1}
     \sum_{i \in \T_n^{\eta}}X_iX_{2i+\eta} - \left(|\T_n^{\eta}|^{-1}
     \sum_{i \in \T_n^{\eta}}X_i\right)\left(|\T_n^{\eta}|^{-1}\sum_{i \in
       \T_n^{\eta}} X_{2i+\eta}\right)}}{\displaystyle{|\T_{n}^{\eta}|^{-1}
     \sum_{i \in \T_{n}^{\eta}}X_i^2 - \left(|\T_{n}^{\eta}|^{-1}\sum_{i \in
       \T_{n}^{\eta}}X_i\right)^2}}, \\ \widehat{\beta}_{\eta}'^{n}  =   \displaystyle{|\T_{n}^{\eta}|^{-1}
   \sum_{i \in \T_{n}^{\eta}} X_{2i+\eta} - \widehat{\alpha}_{\eta}'^{n}
   |\T_{n}^{\eta}|^{-1} \sum_{i \in \T_{n}^{\eta}} X_i}.
  $

Notice that those LSE are based on polynomial functions of the
observations. So, since the latter are bounded, we are in the
functional setting of the results of section \ref{gw_main}.
Recalling the Markov chain $(Y_{n}, n\in\N)$, notice that $Y_{n}$ is
distributed as $Z_{n}=a_{1}a_{2} \cdots a_{n-1}a_{n}  Y_{0}   +
\sum_{k=1}^n   a_{1}a_{2}  \cdots a_{k-1}b_{k}$, where $b_{n}=
b_{n}'  + s_{n} e_{n}$, $((a_{n}, b_{n}', s_{n}), n \ge  1)$ is a
sequence of independent  identically distributed random variables,
whose common distribution is given by, for $\eta \in \{0,1\}$,
\begin{equation*}
\P(a_1=\alpha_{\eta}, b_1'=\beta_{\eta}, s_1=\sigma)=
\frac{p_{1,0}}{m} \quad \mbox{and} \quad \P(a_1=\alpha_{\eta}',
b_1'=\beta_{\eta}', s_1=\sigma_{\eta})= \frac{p_{\eta}}{m},
\end{equation*}
$(e_{n}  ,  n \ge  1)$  is a  sequence  of  independent $\NN(0,1)$
random variables, and is independent of $((a_n, b_n', s_n), n \ge
1)$, and both sequences are independent of $Y_0$. Moreover, it is
easy to check that the sequence $(Z_{n},n\in\N)$ converge a.s. to a
limit $Z$, which implies that the Markov chain $(Y_{n},n\in\N)$
converge in distribution to $Z$. We refer to \cite{DelMar}, section
6, for more details. Following the proof of Proposition 28, step 1
in \cite{Guyon}, we check hypothesis \textbf{(H1)} with $\alpha =
\max(|\alpha_{0}|, |\alpha_{1}|, |\alpha'_{0}|, |\alpha'_{1}|)<1$
and with $\mu$ the distribution of $Z$. Let
$\mu_{1}=\E\left[Z\right]$ and $\mu_{2}=\E\left[Z^{2}\right]$. We
have (see \cite{DelMar})
\begin{equation*}
\mu_{1} = \frac{\overline{\beta}}{1-\overline{\alpha}} \quad
\text{and} \quad \mu_{2} =
\frac{2\overline{\alpha\beta}\overline{\beta}/(1-\overline{\alpha})
+
\overline{\beta^{2}}+\overline{\alpha^{2}}}{1-\overline{\alpha^{2}}},
\end{equation*}
where $\displaystyle \overline{\alpha} = \E\left[a_{1}\right]$,
$\displaystyle \overline{\alpha^{2}} = \E\left[a_{1}^{2}\right]$,
$\displaystyle \overline{\beta} = \E\left[b_{1}\right]$,
$\displaystyle \overline{\beta^{2}} = \E\left[b_{1}^{2}\right]$,
$\displaystyle \overline{\alpha\beta} = \E\left[a_{1}b_{1}\right]$
and $\displaystyle \overline{\sigma^{2}} =
\E\left[s_{1}^{2}\right]$.

We then have the following deviation inequality for
$\widehat{\theta}_{n} - \theta$.

\begin{prop}\label{prop:dev_tgw}
For all $\delta>0$, for all $a>0$, for all $b>0$ and for all
$\gamma>0$ such that $b<a/(\delta+1)$ and
$\gamma<\min\left\{c_{1}/(1+\delta),
c_{1}/\left(1+\sqrt{\delta}\right)\right\}$, where $c_{1}$ is a
positive constant which depends on $p_{1,0}$, $p_{0}$, $p_{1}$,
$\mu_{1}$ and $\mu_{2}$, and for $\displaystyle n_{0} :=
\left(\log\left(\gamma^{q}\delta^{p}b/c_{0}\right)/\log\alpha\right)
- 1$, we have

\begin{itemize}
\item if $m\alpha<1$, then $\forall n\in\N$,
\begin{equation*}
\P\left(\|\widehat{\theta}_{n} - \theta\| > \delta|W\geq
a\right)\leq c_{2}\exp\left(c''\gamma^{q}\delta^{p}b\right)
\exp\left(-c'\left(\gamma^{q}\delta^{p}b\right)^{2}\left(m^{2}/2\right)^{n+1}\right)
+ A_{n};
\end{equation*}
\item if $m\alpha = 1$, then $\forall n\in\N$,
\begin{equation*}
\P\left(\|\widehat{\theta}_{n} - \theta\| > \delta|W\geq
a\right)\leq c_{2}\exp\left(c''\gamma^{q}\delta^{p}b(n+1)\right)
\exp\left(-c'\left(\gamma^{q}\delta^{p}b\right)^{2}\left(m^{2}/2\right)^{n+1}\right)
+ A_{n};
\end{equation*}
\item if $1<m\alpha<\sqrt{2}$, then $\forall n \in\N$ such that $n>n_0$,
\begin{equation*}
\P\left(\|\widehat{\theta}_{n} - \theta\| > \delta|W\geq
a\right)\leq
c_{2}\exp\left(-c'\left(\gamma^{q}\delta^{p}b\right)^{2}\left(m^{2}/2\right)^{n+1}\right)
+ A_{n};
\end{equation*}
\item if $m\alpha=\sqrt{2}$, then $\forall n\in \N$  such that $n>n_0$,
\begin{equation*}
\P\left(\|\widehat{\theta}_{n} - \theta\| > \delta|W\geq
a\right)\leq
c_{2}\exp\left(-c'\left(\gamma^{q}\delta^{p}b\right)^{2}(1/n)\left(m^{2}/2\right)^{n+1}\right)
+ A_{n} ;
\end{equation*}
\item if $m\alpha>\sqrt{2}$, then $\forall n\in \N^{*}$ such that $n>n_0$,
\begin{equation*}
\P\left(\|\widehat{\theta}_{n} - \theta\| > \delta|W\geq
a\right)\leq
c_{2}\exp\left(-c'\left(\gamma^{q}\delta^{p}b\right)^{2}\alpha^{-2n}\right)
+ A_{n};
\end{equation*}
\end{itemize}
where 
$\displaystyle A_{n} =
c_{3}\exp\left(c'\left(\gamma^{q}\delta^{p}b\right)^{2/3}\right)
\exp\left(-c''\left(\gamma^{q}\delta^{p}b\right)^{2/3}\left(t_{n}/(n+1)^{2}\right)^{1/3}\right)$,
$p\in\{1/2, 1\}$, $q\in\{0,1/2,1\}$, $c_{2}$, $c_{3}$, $c_{4}$, $c'$
and $c''$ are positive constants which depend on $c$, $m$, $\alpha$
$p_{1,0}$, $p_{0}$, $p_{1}$, $\mu_{1}$ and $\mu_{2}$.
\end{prop}

\begin{remark}
Note that the constants $c_{2}$, $c_{3}$, $c_{4}$, $c'$ and $c''$
which appear in Proposition \ref{prop:dev_tgw} may differ term by
term. The values of $p$ and $q$ depend on the magnitude of $\delta$
and $\gamma$. For example, for $\delta$ and $\gamma$ small enough,
we have $p=1$ and $q=1$. We also stress that all these constants can
be made explicit by tedious calculations.
\end{remark}

\section{Proofs of the main results}\label{gw_proofs}
\subsection{Proof of Theorem \ref{thm:gw_main1}} Let
$f\in\BB_{b}(S)$ such that $\<\mu,f\> = 0$. We are going to study
successively $\widetilde{M}_{\H_{r}^{*}}(f)$ for $\H_{r}=\G_{r}$ and
$\H_{r}=\T_{r}.$

\textbf{Step 1.} Let us first deal with
$\widetilde{M}_{\G_{r}^{*}}(f)$. By Chernoff inequality, we have for
all $\delta>0$ and for all $\lambda>0$
\begin{equation}\label{cher_g}
\P\left(\widetilde{M}_{\G_{r}^{*}}(f)>\delta\right) \leq
\exp\left(-\lambda\delta m^{r}\right)
\E\left[\exp\left(\lambda\sum\limits_{i\in\G_{r}^{*}}f(X_{i})\right)\right].
\end{equation}
Recall that for all $i\in\G_{r-1}^{*}$,
\begin{equation*}
\E\left[f(X_{2i})\mathbf{1}_{\left\{2i\in\T^{*}\right\}} +
f(X_{2i+1})\mathbf{1}_{\left\{2i+1\in\T^{*}\right\}}|\FF_{r-1}\right]
= mQf(X_{i}).
\end{equation*}
By subtracting and adding terms in expectation of the right hand of
(\ref{cher_g}), and conditioning with respect to $\FF_{r-1},$ we get
\begin{eqnarray}
\E\left[\exp\left(\lambda\sum\limits_{i\in\G_{r}^{*}}f(X_{i})\right)\right]
=
\E\left[\exp\left(\lambda\sum\limits_{i\in\G_{r-1}^{*}}mQf(X_{i})\right)\right.
\hspace{4cm} \label{decomp1}
\\
\times\left.\E\left[\exp\left(\sum\limits_{i\in\G_{r-1}^{*}}\lambda\left(f(X_{2i})\mathbf{1}_{\left\{2i\in\T^{*}\right\}}
+ f(X_{2i+1})\mathbf{1}_{\left\{2i+1\in\T^{*}\right\}} -
mQf(X_{i})\right)\right)\Bigg|\FF_{r-1}\right]\right].\notag
\end{eqnarray}
Observing that $\G_{r-1}^{*}$ is $\FF_{r-1}$ measurable, and using
the fact that conditionally to $\FF_{r-1}$, the triplets
$\{(\Delta_{i}),i\in\G_{r-1}\}$ are independent (this is due to the
Markov property), we have
\begin{align}
\E\left[\exp\left(\sum\limits_{i\in\G_{r-1}^{*}}\lambda\left(f(X_{2i})\mathbf{1}_{\left\{2i\in\T^{*}\right\}}
+ f(X_{2i+1})\mathbf{1}_{\left\{2i+1\in\T^{*}\right\}} -
mQf(X_{i})\right)\right)\Bigg|\FF_{r-1}\right]\hspace{3cm}\label{prod1}\\
=
\prod\limits_{i\in\G_{r-1}^{*}}\E\left[\exp\left(\lambda\left(f(X_{2i})\mathbf{1}_{\left\{2i\in\T^{*}\right\}}
+ f(X_{2i+1})\mathbf{1}_{\left\{2i+1\in\T^{*}\right\}} -
mQf(X_{i})\right)\right)\Bigg|\FF_{r-1}\right].\hspace{2cm}\notag
\end{align}
Using Azuma-Bennet-Hoeffding inequality \cite{Azuma},
\cite{Bennett}, \cite{Hoeffding}, we get according to {\bf (H1)},
for all $i\in\G_{r-1}^{*},$

\vspace{0.25cm}

$\displaystyle
\E\left[\exp\left(\lambda\left(f(X_{2i})\mathbf{1}_{\left\{2i\in\T^{*}\right\}}
+ f(X_{2i+1})\mathbf{1}_{\left\{2i+1\in\T^{*}\right\}} -
mQf(X_{i})\right)\right)\Bigg|\FF_{r-1}\right]\\
\hspace*{8cm} \leq
\exp\left(\frac{c^{2}\lambda^{2}(2+m\alpha)^{2}}{2}\right). $

\vspace{0.25cm}

From (\ref{prod1}), this implies that

\vspace{0.25cm}

$\displaystyle
\E\left[\exp\left(\sum\limits_{i\in\G_{r-1}^{*}}\lambda\left(f(X_{2i})\mathbf{1}_{\left\{2i\in\T^{*}\right\}}
+ f(X_{2i+1})\mathbf{1}_{\left\{2i+1\in\T^{*}\right\}} -
mQf(X_{i})\right)\right)\Bigg|\FF_{r-1}\right]\\ \hspace*{8cm} \leq
\exp\left(\frac{c^{2}\lambda^{2}(2+m\alpha)^{2}|\G_{r-1}^{*}|}{2}\right)
\\ \hspace*{8cm} \leq
\exp\left(\frac{c^{2}\lambda^{2}(2+m\alpha)^{2}|\G_{r-1}|}{2}\right),\hspace{1cm}
$

\vspace{0.25cm}

where we have used the fact that $|\G_{r-1}^{*}|\leq|\G_{r-1}|$ in
the last inequality. Recalling (\ref{decomp1}), we are led to
\begin{eqnarray*}
\E\left[\exp\left(\lambda\sum\limits_{i\in\G_{r}^{*}}f(X_{i})\right)\right]
\leq
\exp\left(\frac{c^{2}\lambda^{2}(2+m\alpha)^{2}|\G_{r-1}|}{2}\right)
\hspace{2cm}
\\ \times
\E\left[\exp\left(\lambda\sum\limits_{i\in\G_{r-1}^{*}}mQf(X_{i})\right)\right].
\end{eqnarray*}
Reproducing the same reasoning with $Qf$ and $\G_{r-1}^{*}$ instead
of $f$ and $\G_{r}^{*},$ we get
\begin{align*}
\E\left[\exp\left(\lambda
m\sum\limits_{i\in\G_{r-1}^{*}}Qf(X_{i})\right)\right] \leq
\exp\left(\frac{c^{2}\lambda^{2}m^{2}(2\alpha+m\alpha^{2})^{2}|\G_{r-2}|}{2}\right)\hspace{4cm}\\
\times\E\left[\exp\left(\lambda
m^{2}\sum\limits_{i\in\G_{r-2}^{*}}Q^{2}f(X_{i})\right)\right].\hspace{2cm}
\end{align*}
Iterating this procedure, we get

\vspace{0.25cm}

$\displaystyle
\E\left[\exp\left(\lambda\sum\limits_{i\in\G_{r}^{*}}f(X_{i})\right)\right]
\leq \exp\left(\frac{c^{2}\lambda^{2}}{2}\sum\limits_{q=0}^{r-1}
\left(2\alpha^{q}+m\alpha^{q+1}\right)^{2}m^{2q}2^{r-1-q}\right)
\\ \hspace*{8cm} \times \E\Big[\exp\left(\lambda
m^{r}Q^{r}f(X_{1})\right)\Big] \\ \hspace*{3cm} \leq
\exp\left(\frac{c^{2}\lambda^{2}(2+m\alpha)^{2}2^{r-1}}{2}\sum\limits_{q=0}^{r-1}
\left(\frac{\alpha^{2}m^{2}}{2}\right)^{q}\right) \times
\exp\left(\lambda c(\alpha m)^{r}\right), $

\vspace{0.25cm}

where the last inequality was obtained from {\bf (H1)}. From the
foregoing and from (\ref{cher_g}), we deduce that

\vspace{0.25cm}
$\displaystyle \P\left(\widetilde{M}_{\G_{r}^{*}}(f)>\delta\right)
\leq
\begin{cases}
\exp\left(-\lambda\delta m^{r} +
\frac{c^{2}\lambda^{2}(2+m\alpha)^{2}\left(2^{r}-(\alpha^{2}m^{2})^{r}\right)}{2(2-\alpha^{2}m^{2})}\right)
\\ \hspace*{3.5cm} \times \exp\left(\lambda c(\alpha m)^{r}\right) \hspace{0.25cm} \text{if
$\alpha^{2}m^{2} \neq 2$}, \\  \exp\left(-\lambda\delta m^{r} +
c^{2}\lambda^{2}(2+\sqrt{2})^{2}r2^{r-2}\right) \exp\left(\lambda
c(\sqrt{2})^{r}\right) \\ \hspace*{7cm} \text{if $\alpha^{2}m^{2}
=2$}.
\end{cases}
$

\vspace{0.25cm}

Now, the rest divides into four cases. In the sequel $c_{1}$ and
$c_{2}$ will denote positive constants which depend on $c$, $m$, and
$\alpha$.

$\bullet$ If $m\alpha\leq 1$, then, for all $r\in\N$,
$(m\alpha)^{r}<1$ and $2^{r}-(\alpha^{2}m^{2})^{r}<2^{r}$. We then
have
\begin{equation*}
\P\left(\widetilde{M}_{\G_{r}^{*}}(f)>\delta\right) \leq
\exp\left(c\lambda\right)\exp\left(-\lambda\delta m^{r} +
\lambda^{2} c_{1} 2^{r}\right).
\end{equation*}
Taking $\lambda = (\delta m^{r})/(2^{r+1}c_{1})$, we are led to
\begin{equation*}
\P\left(\widetilde{M}_{\G_{r}^{*}}(f)>\delta\right) \leq
\exp\left(c_{1}\delta\right)\exp\left(-\delta^{2}c_{1}\left(\frac{m^{2}}{2}\right)^{r}\right).
\end{equation*}

$\bullet$ If $1<m\alpha<\sqrt{2},$ then, since
$2^{r}-(\alpha^{2}m^{2})^{r}<2^{r}$, we have
\begin{equation*}
\P\left(\widetilde{M}_{\G_{r}^{*}}(f)>\delta\right) \leq
\exp\left(-\lambda\delta m^{r} +
\lambda^{2}c_{1}2^{r}\right)\exp\left(\lambda c(m\alpha)^{r}\right).
\end{equation*}
Taking $\lambda = (\delta m^{r})/(2^{r+1}c_{1})$, we are led to
\begin{equation*}
\P\left(\widetilde{M}_{\G_{r}^{*}}(f)>\delta\right) \leq
\exp\left(-c_{2}\delta (m^{2}/2)^{r}(\delta - 2c\alpha^{r})\right).
\end{equation*}
For all $r\in\N$ such that $r>\log(\delta/4c)/\log(\alpha),$ we have
$\delta-2c\alpha^{r}>\delta/2$ and it then follows that
\begin{equation*}
\P\left(\widetilde{M}_{\G_{r}^{*}}(f)>\delta\right) \leq
\exp\left(-c_{2}\delta^{2}(m^{2}/2)^{r}\right).
\end{equation*}
$\bullet$ If $m\alpha=\sqrt{2}$, then we have
\begin{equation*}
\P\left(\widetilde{M}_{\G_{r}^{*}}(f)>\delta\right) \leq
\exp\left(-\lambda\delta
m^{r}+\lambda^{2}c_{1}r2^{r-2})\right)\exp\left(\lambda
c\left(\sqrt{2}\right)^{r}\right).
\end{equation*}
Taking $\lambda=(\delta m^{r})/(c_{1}r2^{r-1})$, we have for all
$r>\log(\delta/4c)/\log(\sqrt{2}/m)$,
\begin{equation*}
\P\left(\widetilde{M}_{\G_{r}^{*}}(f)>\delta\right) \leq
\exp\left(-c_{2}\delta^{2}(1/r)(m^{2}/2)^{r}\right).
\end{equation*}
$\bullet$ If $m\alpha>\sqrt{2}$, then we have
\begin{equation*}
\P\left(\widetilde{M}_{\G_{r}^{*}}(f)>\delta\right) \leq
\exp\left(-\lambda\delta m^{r} + \lambda^{2}
c_{1}(m^{2}\alpha^{2})^{r}\right)\exp\left(\lambda
c(m\alpha)^{r}\right).
\end{equation*}
Taking $\lambda = \delta/(2c_{1}(m\alpha^{2})^{r}),$ we have for all
$r>\log(\delta/4c)/\log\alpha$,
\begin{equation*}
\P\left(\widetilde{M}_{\G_{r}^{*}}(f)>\delta\right) \leq
\exp\left(-c_{3}\delta^{2}\alpha^{-2r}\right).
\end{equation*}
This ends the proof for $\H_{r}=\G_{r}$.


\textbf{Step 2.} Let us look at $\widetilde{M}_{\T_{r}^{*}}(f)$. By
Chernoff inequality, we have for all $\delta>0$ and for all
$\lambda>0$
\begin{equation}\label{cher_t}
\P\left(\widetilde{M}_{\T_{r}^{*}}(f)>\delta\right) \leq
\exp\left(-\lambda\delta
t_{r}\right)\E\left[\exp\left(\lambda\sum\limits_{i\in\T_{r}^{*}}f(X_{i})\right)\right].
\end{equation}
Expectation which appears in the right hand of (\ref{cher_t}) can be
written as
\begin{eqnarray}\label{decomp2}
\E\left[\exp\left(\lambda\sum\limits_{i\in\T_{r}^{*}}
f(X_{i})\right)\right] =
\E\left[\exp\left(\lambda\sum\limits_{i\in\T_{r-2}^{*}}
f(X_{i})\right)\right. \hspace{5cm} \\ \times
\exp\left(\lambda\sum\limits_{i\in\G_{r-1}^{*}}
(f+mQf)(X_{i})\right) \notag  \\
\times \left.
\E\left[\exp\left(\lambda\sum\limits_{i\in\G_{r-1}^{*}}
\left(f(X_{2i})\mathbf{1}_{\{2i\in\T^{*}\}} +
f(X_{2i+1})\mathbf{1}_{\{2i+1\in\T^{*}\}} -
mQf(X_{i})\right)\right)\bigg|\FF_{r-1}\right]\right].\notag
\end{eqnarray}
Observing that $\G_{r-1}^{*}$ is $\FF_{r-1}$ measurable, and using
the fact that conditionally to $\FF_{r-1}$, the triplets
$\{(\Delta_{i}),i\in\G_{r-1}\}$ are independent and
Azuma-Bennet-Hoeffding inequality, we obtain

\vspace{0.25cm}

$\displaystyle
\E\left[\exp\left(\lambda\sum\limits_{i\in\G_{r-1}^{*}}
\left(f(X_{2i})\mathbf{1}_{\{2i\in\T^{*}\}} +
f(X_{2i+1})\mathbf{1}_{\{2i+1\in\T^{*}\}} -
mQf(X_{i})\right)\right)\bigg|\FF_{r-1}\right]  \\ \hspace*{0.75cm}
= \prod\limits_{i\G_{r-1}^{*}} \E\left[\exp\left(\lambda
\left(f(X_{2i})\mathbf{1}_{\{2i\in\T^{*}\}} +
f(X_{2i+1})\mathbf{1}_{\{2i+1\in\T^{*}\}} -
mQf(X_{i})\right)\right)\bigg|\FF_{r-1}\right] \\ \hspace*{2cm} \leq
\exp\left(\frac{c^{2}\lambda^{2}(2+m\alpha)^{2}|\G_{r-1}^{*}|}{2}\right)
\\ \hspace*{2cm} \leq \exp\left(\frac{c^{2}(2+m\alpha)^{2}|\G_{r-1}|}{2}\right), \hspace{7.5cm}
$

\vspace{0.25cm}

where the last inequality was obtained using the fact that
$|\G_{r-1}^{*}|\leq |\G_{r-1}|$. From the foregoing and from
(\ref{decomp2}), we deduce that

\vspace{0.25cm}

$\displaystyle \E\left[\exp\left(\lambda\sum\limits_{i\in\T_{r}^{*}}
f(X_{i})\right)\right] \leq
\exp\left(\frac{c^{2}(2+m\alpha)^{2}|\G_{r-1}|}{2}\right) \\
\hspace*{3cm} \times
\E\left[\exp\left(\lambda\sum\limits_{i\in\T_{r-2}^{*}}
f(X_{i})\right) \exp\left(\lambda\sum\limits_{i\in\G_{r-1}^{*}}
(f+mQf)(X_{i})\right)\right] $

\vspace{0.25cm}

Doing the same thing with $(f+mQf)$ and $\G_{r-1}^{*}$ instead of
$f$ and $\G_{r}^{*},$ we get

\vspace{0.25cm}

$\displaystyle
\E\left[\exp\left(\lambda\sum\limits_{i\in\T_{r-2}^{*}}
f(X_{i})\right) \exp\left(\lambda\sum\limits_{i\in\G_{r-1}^{*}}
(f+mQf)(X_{i})\right)\right] \\ \hspace*{1cm} =
\E\left[\exp\left(\lambda\sum\limits_{i\in\T_{r-3}^{*}}
f(X_{i})\right) \times
\exp\left(\lambda\sum\limits_{i\in\G_{r-2}^{*}}
(f+mQf+m^{2}Q^{2}f)(X_{i})\right)\right.  \\ \hspace*{0.75cm} \times
\E\left[\exp\left(\lambda\sum\limits_{i\in\G_{r-2}^{*}}
\bigg((f+mQf)(X_{2i})\mathbf{1}_{\{2i\in\T^{*}\}} +
(f+mQf)(X_{2i+1})\mathbf{1}_{\{2i+1\in\T^{*}\}} \bigg. \right.
\right. \\ \hspace*{7.5cm} - \Bigg. \Bigg. \Bigg. \bigg.
(mQf+m^{2}Q^{2}f)(X_{i})\bigg)\Bigg)\bigg|\FF_{r-1}\Bigg] \Bigg]
\\ \hspace*{1cm} \leq
\exp\left(\frac{c^{2}\lambda^{2}(2+3m\alpha+m^{2}\alpha^{2})^{2}|\G_{r-2}|}{2}\right)
\\ \hspace*{1.5cm} \times \E\left[\exp\left(\lambda\sum\limits_{i\in\T_{r-3}^{*}}
f(X_{i})\right) \times
\exp\left(\lambda\sum\limits_{i\in\G_{r-2}^{*}}
(f+mQf+m^{2}Q^{2}f)(X_{i})\right)\right]. $

\vspace{0.25cm}

Iterating this procedure, we are led to

\vspace{0.25cm}

$\displaystyle \E\left[\exp\left(\lambda\sum\limits_{i\in\T_{r}^{*}}
f(X_{i})\right)\right] \leq
\exp\left(\frac{c^{2}(2+m\alpha)^{2}\lambda^{2}}{2}
\sum\limits_{q=1}^{r} \left(\sum\limits_{k=0}^{q-1}
(m\alpha)^{k}\right)^{2}2^{r-q}\right) \\ \hspace*{8cm} \times
\E\left[\exp\left(\lambda\sum\limits_{q=0}^{r}
m^{q}Q^{q}f(X_{1})\right)\right] \hspace{2cm} \\ \hspace*{2cm} \leq
\exp\left(\frac{c^{2}(2+m\alpha)^{2}\lambda^{2}}{2}
\sum\limits_{q=1}^{r} \left(\sum\limits_{k=0}^{q-1}
(m\alpha)^{k}\right)^{2}2^{r-q}\right) \exp\left(\lambda
c\sum\limits_{q=0}^{r} (m\alpha)^{q}\right), $

\vspace{0.25cm}

where the last inequality was obtained using hypothesis
\textbf{(H1)}. In the sequel, $c_{0}$, $c_{1}$ and $c_{2}$ will
denote some positive constants which depend on $\alpha$, $m$, and
$c$. They may differ from one line to another. For $m\alpha\neq 1$
and $m\alpha\neq \sqrt{2}$, we deduce from the foregoing and from
(\ref{cher_t}) that

\vspace{0.25cm}

$\displaystyle \P\left(\widetilde{M}_{\T_{r}^{*}}(f)>\delta\right)
\leq \exp\left(-\lambda\delta t_{r}\right) \exp\left(\frac{\lambda c
(1-(m\alpha)^{r+1})}{1-m\alpha}\right) \\ \hspace*{1cm} \times
\exp\left(\frac{c^{2}(2+m\alpha)^{2}\lambda^{2}}{2(1-m\alpha)^{2}}\left((2^{r}-1)
- \frac{2m\alpha(2^{r}-(m\alpha)^{r})}{2-m\alpha}\right.\right.
\\ \hspace*{9cm} \left.\left. + \frac{(m\alpha)^{2}
(2^{r}-(m^{2}\alpha^{2})^{r})}{2-(m\alpha)^{2}}\right)\right)
\\ \hspace*{1cm} \leq \exp\left(-\lambda\delta t_{r}+
\frac{c^{2}(2+m\alpha)^{2}\lambda^{2}}{2(m\alpha-1)^{2}}
\left((2^{r}-1) +
\frac{(m\alpha)^{2}(2^{r}-(m^{2}\alpha^{2})^{r})}{2-(m\alpha)^{2}}\right)\right)
\\ \hspace*{8.5cm} \times \exp\left(\frac{\lambda c
(1-(m\alpha)^{r+1})}{1-m\alpha}\right).$

\vspace{0.25cm}

Taking $\displaystyle \lambda = \frac{\delta
t_{r}(m\alpha-1)^{2}}{c^{2}(2+m\alpha)^{2}\left((2^{r}-1) +
\frac{(m\alpha)^{2}(2^{r}-(m^{2}\alpha^{2})^{r})}{2-(m\alpha)^{2}}\right)}$,
we are led to
\vspace{0.25cm}

$\displaystyle \P\left(\widetilde{M}_{\T_{r}^{*}}(f)>\delta\right)
\leq
\exp\left(-\frac{\delta^{2}(1-m\alpha)^{2}t_{r}^{2}}{2c^{2}(2+m\alpha)^{2}\left(2^{r}-1
+
\frac{(m\alpha)^{2}(2^{r}-(m^{2}\alpha^{2})^{r})}{2-(m\alpha^){2}}\right)}\right)
\\ \hspace*{1cm} \times \exp\left(\frac{\delta(1-m\alpha)^{2}t_{r}}{c(2+m\alpha)^{2}\left(2^{r}-1 +
\frac{(m\alpha)^{2}(2^{r}-(m^{2}\alpha^{2})^{r})}{2-(m\alpha)^{2}}\right)}
\times \frac{1-(m\alpha)^{r+1}}{1-m\alpha}\right). $
\vspace{0.25cm}

Now, the rest of the proof divides into five cases.

$\bullet$ If $m\alpha<1$, then, for all $r\in\N$,
$(m\alpha)^{r+1}-1\leq m\alpha -1$ and $2^{r}-(m\alpha)^{2r}<2^{r}$.
We then deduce that
\begin{equation*}
\P\left(\widetilde{M}_{\T_{r}^{*}}(f)>\delta\right) \leq
\exp\left(c_{2}\delta\right)
\exp\left(-c_{2}\delta^{2}(m^{2}/2)^{r+1}\right).
\end{equation*}

$\bullet$ If $1<m\alpha<\sqrt{2}$, then we have
\vspace{0.25cm}

$\displaystyle \hspace*{2cm}
\P\left(\widetilde{M}_{\T_{r}^{*}}(f)>\delta\right) \leq
\exp\left(-c_{1}\delta^{2}(m^{2}/2)^{r+1}\right)\exp\left(c_{2}\delta
\frac{(m\alpha)^{r+1}-1}{m\alpha-1}\right) \\ \hspace*{4.5cm} \leq
\exp\left(-\delta
c_{2}(m^{2}/2)^{r+1}(\delta-c_{0}\alpha^{r+1})\right)$.

\vspace{0.25cm}
Now, for all $r\in\N$ such that
$r+1>\log(\delta/2c_{0})/\log(\alpha),$ we have
$\delta-c_{0}\alpha^{r+1}>\delta/2$, in such a way that
\begin{equation*}
\P\left(\widetilde{M}_{\T_{r}^{*}}(f)>\delta\right) \leq
\exp\left(\delta^{2}c_{2}(m^{2}/2)^{r+1}\right).
\end{equation*}

$\bullet$ If $m\alpha>\sqrt{2}$, then for all $r\in\N$,
$(m^{2}\alpha^{2})^{r}>2^{r}$. We then have
\begin{equation*}
\P\left(\widetilde{M}_{\T_{r}^{*}}(f)>\delta\right) \leq
\exp\left(-c_{2}\delta\alpha^{-2r}(\delta-c_{0}\alpha^{r+1})\right).
\end{equation*}
Now for all $r\in\N$ such that
$r+1>\log(\delta/c_{0})/\log(\alpha),$ we have
\begin{equation*}
\P\left(\widetilde{M}_{\T_{r}^{*}}(f)>\delta\right) \leq
\exp\left(-\frac{c_{2}\delta^{2}}{\alpha^{2r}}\right).
\end{equation*}

$\bullet$ If $m\alpha=1$, then
\vspace{0.25cm}

$\displaystyle \P\left(\widetilde{M}_{\T_{r}^{*}}(f)>\delta\right)
\leq \exp\left(-\lambda\delta t_{r} +
c_{1}2^{r}\lambda^{2}\right)\exp\left(\lambda c(r+1)\right)$
\vspace{0.25cm}

Taking $\displaystyle\lambda = \delta t_{r}/c_{1}2^{r+1}$, we are
led to
\vspace{0.25cm}

$\displaystyle \hspace*{2cm}
\P\left(\widetilde{M}_{\T_{r}^{*}}(f)>\delta\right) \leq
\exp\left(c_{1}\delta \frac{(r+1)t_{r}}{2^{r+1}}\right)
\exp\left(-c_{2}\delta^{2}(m^{2}/2)^{r+1}\right). $

\vspace{0.25cm}

$\bullet$ If $m\alpha=\sqrt{2}$, then
\vspace{0.25cm}

$\displaystyle
\hspace*{2cm}\P\left(\widetilde{M}_{\T_{r}^{*}}(f)>\delta\right)
\leq \exp\left(-\lambda\delta t_{r} + \lambda^{2}
c_{1}(r+1)2^{r}\right)\exp\left(\lambda
c_{1}(\sqrt{2})^{r+1}\right)$.

\vspace{0.25cm}
Taking $\displaystyle \lambda=\delta t_{r}/(2c_{1}(r+1)2^{r})$, we
are led to
\vspace{0.25cm}

$\displaystyle
\hspace*{2cm}\P\left(\widetilde{M}_{\T_{r}^{*}}(f)>\delta\right)
\leq \exp\left(-\frac{c_{2}\delta}{r+1}
\left(\frac{m^{2}}{2}\right)^{r+1}
\left(\delta-c_{0}\left(\frac{\sqrt{2}}{m}\right)^{r+1}\right)\right).
$

\vspace{0.25cm}
Now, for all $r\in\N$ such that
$r+1>\log(\delta/c_{0})/\log(\sqrt{2}/m)$, we get
\vspace{0.25cm}

$\displaystyle
\hspace*{2cm}\P\left(\widetilde{M}_{\T_{r}^{*}}(f)>\delta\right)
\leq
\exp\left(-\frac{c_{2}\delta^{2}}{r+1}\left(\frac{m^{2}}{2}\right)^{r+1}\right)
$.

\vspace{0.25cm}
This ends the proof for $\H_{r}=\T_{r}$.

%
%
%
%

\subsection{Proof of Theorem \ref{thm:gw_main11}} Let
$f\in\mathcal{B}_{b}(S)$ such that $\<\mu,f\>\neq0$. Once again, we
are going to study successively $\widetilde{M}_{\G_{r}^{*}}(f)$ and
$\widetilde{M}_{\T_{r}^{*}}(f)$.

\textbf{Step 1.} Let us first deal with
$\widetilde{M}_{\G_{r}^{*}}(f)$. Set $g=f-\<\mu,f\>.$ Then,
$\<\mu,g\>=0$ and
\begin{equation*}
\widetilde{M}_{\G_{r}^{*}}(f) = \widetilde{M}_{\G_{r}^{*}}(g) +
(|\G_{r}^{*}|/m^{r})\<\mu,f\>.
\end{equation*}
We have
\begin{eqnarray}\label{ineq1}
\P\left(\widetilde{M}_{\G_{r}^{*}}(f)-\<\mu,f\>W>\delta\right) \leq
\P\bigg(\widetilde{M}_{\G_{r}^{*}}(g)>\delta/2\bigg) \hspace{3cm} \\
+
\P\left(\left|\frac{|\G_{r}^{*}|}{m^{r}}-W\right|>\frac{\delta}{2|\<\mu,f\>|}\right).\notag
\end{eqnarray}
As $\<\mu,g\>=0,$ the previous computations (proof of Theorem
\ref{thm:gw_main1}) give us some bound for the first term of right
hand of (\ref{ineq1}), similar to those obtain in Theorem
\ref{thm:gw_main1}. Now, under hypothesis {\bf (H3)}, we deduce,
from \cite{At94} Theorem 5, that
\begin{equation*}
\P\left(\left|\frac{|\G_{r}^{*}|}{m^{r}}-W\right|>\frac{\delta}{2|\<\mu,f\>|}\right)
\leq c_{2}\exp\left(-c_{3}\delta^{2/3}m^{r/3}\right),
\end{equation*}
and this ends the proof of Theorem \ref{thm:gw_main11} when
$\H_{r}=\G_{r}$.

\textbf{Step 2.} Let us look at $\widetilde{M}_{\T_{r}^{*}}(f)$. For
$f\in\BB_{b}(S),$ set $g=f-\<\mu,f\>$. Then, $\<\mu,g\>=0$ and
\begin{equation*}
\widetilde{M}_{\T_{r}^{*}}(f) = \widetilde{M}_{\T_{r}^{*}}(g) +
(|\T_{r}^{*}|/t_{r})\<\mu,f\>.
\end{equation*}
We have
\begin{eqnarray}\label{ineq2}
\P\left(\widetilde{M}_{\T_{r}^{*}}(f)-\<\mu,f\>W>\delta\right) \leq
\P\left(\widetilde{M}_{\T_{r}^{*}}(g)>\delta/2\right) \hspace{3cm} \\
+ \P\left(\left|\frac{|\T_{r}^{*}|}{t_{r}}-W\right| >
\frac{\delta}{2|\<\mu,f\>|}\right).\notag
\end{eqnarray}
Since $\<\mu,g\>=0,$ the first term of the right hand of
(\ref{ineq2}) can be bounded as in the previous computations (proof
of Theorem \ref{thm:gw_main1}). Under additional hypothesis
\textbf{(H3)}, we have, from \cite{At94} Theorem 5,

\vspace{0.25cm}

$\displaystyle \P\left(\left|\frac{|\T_{r}^{*}|}{t_{r}}-W\right|
> \frac{\delta}{2|\<\mu,f\>|}\right) \leq \sum\limits_{q=0}^{r}
\P\left(\frac{m^{q}}{t_{r}}\left|\frac{|\G_{q}^{*}|}{m^{q}}-W\right|
> \frac{\delta}{2(r+1)|\<\mu,f\>|}\right) \\ \hspace*{6cm} =\sum\limits_{q=0}^{r}
\P\left(\left|\frac{|\G_{q}^{*}|}{m^{q}}-W\right|
> \frac{\delta t_{r}}{2(r+1)|\<\mu,f\>|m^{q}}\right) \\ \hspace*{6cm} \leq
\sum\limits_{q=0}^{r}c_{2}
\exp\left(-c_{3}\delta^{2/3}\left(\frac{t_{r}^{2}}{(r+1)m^{q}}\right)^{1/3}\right)
\\ \hspace*{6cm} \leq
c_{2}\exp\left(-c_{3}\delta^{2/3}\left(\frac{t_{r}}{(r+1)^{2}}\right)^{1/3}\right)\bigg(1+o(1)\bigg),
$

\vspace{0.25cm}

and this ends the proof of Theorem \ref{thm:gw_main11} when
$\H_{r}=\T_{r}.$

\subsection{Proof of Theorem \ref{thm:gw_main2}}
Let $f\in\BB_{b}(S)$. Without loss of generality, we assume that
$\<\mu,f\>=0$. Otherwise, we take $f-\langle\mu,f\rangle$. For all
$\delta>0$, for all $a>0$ and for all $b>0$ such that
$b<a/(\delta+1),$ we have
\vspace{0.25cm}

$\displaystyle \P\left(\overline{M}_{\H_{r}^{*}}(f)>\delta|W\geq
a\right) = \P\left(\overline{M}_{\H_{r}^{*}}(f)>\delta,
\frac{|\H_{r}^{*}|}{h_{r}}>b|W\geq a\right) \\ \hspace*{5cm} +
\P\left(\overline{M}_{\H_{r}^{*}}(f)>\delta,
\frac{|\H_{r}^{*}|}{h_{r}}\leq b|W\geq a\right) \\ \hspace*{4cm} =
\frac{1}{\P\left(W\geq
a\right)}\Bigg(\P\left(\overline{M}_{\H_{r}^{*}}(f)>\delta,
\frac{|\H_{r}^{*}|}{h_{r}}>b, W\geq a\right) \Bigg.  \\
\hspace*{5cm} + \Bigg. \P\left(\overline{M}_{\H_{r}^{*}}(f)>\delta,
\frac{|\H_{r}^{*}|}{h_{r}}\leq b, W\geq a\right)\Bigg) \\
\hspace*{2cm} \leq p_{a}
\P\left(\widetilde{M}_{\H_{r}^{*}}(f)>\delta b\right) + p_{a}
\P\left(\left|\frac{|\H_{r}^{*}|}{h_{r}}-W\right|>W-b,W\geq a\right)
\\ \\ \hspace*{2cm} \leq p_{a} \P\left(\widetilde{M}_{\H_{r}^{*}}(f)>\delta b\right)
+ p_{a} \P\left(\left|\frac{|\H_{r}^{*}|}{h_{r}}-W\right|> \delta
b\right)$,

\vspace{0.25cm}
where $p_{a}=\P\left(W\geq a\right)^{-1}$. Now, the first term of
the last inequality can be bounded as in  Theorem
\ref{thm:gw_main1}, and the second term is bounded as in the
\textbf{step 1} and and \textbf{step 2} of the proof of Theorem
\ref{thm:gw_main11}. This ends the proof.

\subsection{Proof of Theorem \ref{thm:gw_main3}}
Let $f\in\BB_{b}(S^{3})$.

\textbf{Step 1.} Let us first deal with
$\widetilde{M}_{\G_{r}^{*}}(f).$ Assume that $\<\mu,P^{*}f\>=0.$ By
Chernoff inequality, we have for all $\delta>0$ and for all
$\lambda>0,$
\begin{equation*}
\P\left(\widetilde{M}_{\G_{r}^{*}}(f)>\delta\right) \leq
\exp\left(-\lambda\delta m^{r}\right) \E\left[\exp\left(\lambda
\sum\limits_{i\in\G_{r}^{*}} f(\Delta_{i})\right)\right].
\end{equation*}
Conditioning by $\FF_{r},$ and using, conditional independence of
triplets $\{\Delta_{i},i\in\G_{r}\}$ with respect to $\FF_{r}$,
Azuma-Bennet-Hoeffding inequality and \textbf{(H2)}, we obtain

\vspace{0.25cm}

$\displaystyle \E\left[\exp\left(\lambda
\sum\limits_{i\in\G_{r}^{*}} f(\Delta_{i})\right)\right] \\
\hspace*{1cm} = \E\left[\exp\left(\lambda\sum\limits_{i\G_{r}^{*}}
P^{*}f(X_{i})\right)
\E\left[\exp\left(\lambda\sum\limits_{i\in\G_{r}^{*}} (f(\Delta_{i})
- P^{*}f(X_{i}))\right)\Big|\FF_{r}\right]\right] \\ \hspace*{1cm} =
\E\left[\exp\left(\lambda\sum\limits_{i\G_{r}^{*}}
P^{*}f(X_{i})\right)
\prod\limits_{i\in\G_{r}^{*}}\E\left[\exp\left(\lambda(f(\Delta_{i})
- P^{*}f(X_{i}))\right)\Big|\FF_{r}\right]\right] \\ \hspace*{1cm}
\leq \exp\left(2\lambda^{2}\|f\|_{\infty}c_{1}m^{r}\right)
\E\left[\exp\left(\lambda\sum\limits_{i\in\G_{r}^{*}}
P^{*}f(X_{i})\right)\right].$

\vspace{0.25cm}

We control the last expectation as in the \textbf{Step 1} of the
proof of Theorem \ref{thm:gw_main1}, apply to $P^{*}f.$ Next, we get
the result discussing as in the proof of Theorem \ref{thm:gw_main1}.

If $\<\mu,P^{*}f\>\neq 0,$ we set $g=f-\<\mu,P^{*}f\>.$ Then, we
have
\begin{eqnarray}\label{ineq3}
\hspace*{-0.15cm} \P\left(\widetilde{M}_{\G_{r}^{*}}(f) -
\<\mu,P^{*}f\>W>\delta\right) \leq
\P\left(\widetilde{M}_{\G_{r}^{*}}(g)>\delta/2\right) \hspace{3cm} \\
+
\P\left(\left|\frac{\G_{r}^{*}}{m^{r}}-W\right|>\delta/2|\<\mu,P^{*}f\>|\right).\notag
\end{eqnarray}
The first term of the right hand of (\ref{ineq3}) can be bounded as
previously since $\<\mu,P^{*}g\>=0.$ The second term can be bounded
as in \textbf{Step 1} of the proof of Theorem \ref{thm:gw_main1}.
This ends the proof for $\widetilde{M}_{\G_{r}^{*}}(f).$

\textbf{Step 2.} Let us now treat $\widetilde{M}_{\T_{r}^{*}}(f).$
First, we assume that $\<\mu,P^{*}f\>=0.$ For all $\delta>0,$ we
have
\begin{equation*}
\P\left(\widetilde{M}_{\T_{r}^{*}}(f)>\delta\right) \leq
\P\left(\frac{1}{t_{r}}\sum\limits_{i\in\T_{r}^{*}}
(f(\Delta_{i})-P^{*}f(X_{i})) >\delta/2\right) +
\P\left(\widetilde{M}_{\T_{r}^{*}}(P^{*}f)>\delta/2\right).
\end{equation*}
By chernoff inequality, we have for all $\lambda>0,$
\begin{eqnarray*}
\P\left(\frac{1}{t_{r}}\sum\limits_{i\in\T_{r}^{*}}
(f(\Delta_{i})-P^{*}f(X_{i})) > \delta/2\right) \leq
\exp\left(-\frac{\lambda\delta t_{r}}{2}\right) \hspace{3cm} \\
\times \E\left[\exp\left(\lambda\sum\limits_{i\in\T_{r}^{*}}
(f(\Delta_{i})-P^{*}f(X_{i}))\right)\right]
\end{eqnarray*}
Conditioning successively with respect to $(\FF_{q})_{0\leq q\leq
r},$ using conditional independence of triplets
$\{\Delta_{i},i\in\G_{q}\}$ with respect to $\FF_{q}$ and applying
successively Azuma-Bennet-Hoeffding inequality and the fact that
$|\G_{q}^{*}|\leq|\G_{q}|$ for all $q\in\{0,\cdots, r\}$, we get

\vspace{0.25cm}

$\displaystyle \E\left[\exp\left(\lambda\sum\limits_{i\in\T_{r}^{*}}
(f(\Delta_{i})-P^{*}f(X_{i}))\right)\right] \\ \hspace*{.5cm} =
\E\left[\exp\left(\lambda\sum\limits_{i\in\T_{r-1}^{*}}
(f(\Delta_{i})-P^{*}f(X_{i}))\right)\E\left[\exp\left(\lambda\sum\limits_{i\in\G_{r}^{*}}
(f(\Delta_{i})-P^{*}f(X_{i}))\right)\Big|\FF_{r}\right]\right] \\
\hspace*{.5cm} =
\E\left[\exp\left(\lambda\sum\limits_{i\in\T_{r-1}^{*}}
(f(\Delta_{i})-P^{*}f(X_{i}))\right)\prod\limits_{i\in\G_{r}^{*}}\E\left[\exp\left(\lambda
(f(\Delta_{i})-P^{*}f(X_{i}))\right)\Big|\FF_{r}\right]\right] \\
\hspace*{.5cm} \leq
\exp\left(2\lambda^{2}\|f\|_{\infty}^{2}|\G_{r}|\right)
\E\left[\exp\left(\lambda\sum\limits_{i\in\T_{r-1}^{*}}
(f(\Delta_{i})-P^{*}f(X_{i}))\right)\right] \\ \hspace*{3cm} \vdots \\
\hspace*{.5cm} \leq
\exp\left(2\lambda^{2}\|f\|_{\infty}^{2}|\T_{r}|\right).$

\vspace{0.25cm}

Next, optimizing on $\lambda$, we obtain
\begin{equation*}
\P\left(\frac{1}{t_{r}}\sum\limits_{i\in\T_{r}^{*}}
(f(\Delta_{i})-P^{*}f(X_{i})) >\delta/2\right) \leq
\exp\left(-c_{1}\delta^{2}\left(\frac{m^{2}}{2}\right)^{r+1}\right),
\end{equation*}
for some positive constant $c_{1}$. The term
$\P\left(\widetilde{M}_{\T_{r}^{*}}(P^{*}f)>\delta/2\right)$ can be
bounded as in the proof of Theorem \ref{thm:gw_main1}, and this ends
the proof when $\<\mu,P^{*}f\>=0.$ On the other hand, if
$\<\mu,P^{*}f\>\neq 0,$ we have
\begin{equation*}
\widetilde{M}_{\T_{r}^{*}}(f) - \<\mu,P^{*}f\>W =
\widetilde{M}_{\T_{r}^{*}}(g) + \left(\frac{|\T_{r}^{*}|}{t_{r}} -
W\right)\<\mu,P^{*}f\>.
\end{equation*}
We then proceed as for (\ref{ineq3}), and this ends the proof for
$\widetilde{M}_{\T_{r}^{*}}(f).$

\textbf{Step 3.} Eventually, we bound
$\P\left(\overline{M}_{\H_{r}^{*}}(f)>\delta-\<\mu,P^{*}f\>
>\delta\right),$ using \textbf{Step 1} and \textbf{Step 2}, as in the proof of Theorem \ref{thm:gw_main2}.

\subsection{Proof of Proposition \ref{prop:dev_tgw}}
We are going to treat $\widehat{\alpha}_{0}^{n}-\alpha_{0}$.
Deviation inequalities for $\widehat{\alpha}_{1}^{n}-\alpha_{1}$,
$\widehat{\beta}_{\eta}^{n}-\beta_{\eta}$,
$\widehat{\alpha}_{\eta}^{'n}-\alpha_{\eta}^{'}$,
$\widehat{\beta}_{\eta}^{'n}-\beta_{\eta}^{'}$, $\eta\in\{0,1\}$,
can be treated in the same way. Recalling that the state space of
the process $(X_{i},i\in\T^{*})$, denoted by $S$, is assumed to be a
compact subset of $\R$.

Let $g_{1}$, $g_{2}$, $h_{1}$ and $h_{2}$ the functions defined on
$S^{3}$ respectively by
\begin{center}
$\displaystyle g_{1}(x,y,z) = \left(xy - x(\alpha_{0}x +
\beta_{0})\right)\mathbf{1}_{S^{3}}(x,y,z)$,
\end{center}
\begin{center}
$\displaystyle g_{2}(x,y,z) = \left(y - \alpha_{0}x -
\beta_{0})\right)\mathbf{1}_{S^{3}}(x,y,z)$,
\end{center}
\begin{center}
$\displaystyle h_{1}(x,y,z) = x\mathbf{1}_{S^{3}}(x,y,z)$,
\end{center}
\begin{center}
$\displaystyle h_{2}(x,y,z) = x^{2}\mathbf{1}_{S^{3}}(x,y,z)$.
\end{center}
It is easy to see that $\displaystyle P^{*}g_{1}(x) = 0$,
$\displaystyle P^{*}g_{2}(x) = 0$, $\displaystyle P^{*}h_{1}(x) =
p_{1,0}x$ and $\displaystyle P^{*}h_{2}(x) = p_{1,0}x^{2}$ where
$P^{*}$ denote the transition kernel associated to the BAR(1)
process with missing data. With these notations, we can rewrite
$\widehat{\alpha}_{0}^{n}-\alpha_{0}$ as
\vspace{0.25cm}

$\displaystyle \widehat{\alpha}_{0}^{n}-\alpha_{0} =
\frac{|\T_{n}^{*}|^{-1}|\T_{n}^{1,0}|\left(
|\T_{n}^{*}|^{-1}\sum_{i\in\T_{n}^{*}}
g_{1}(\Delta_{i})\right)}{B_{n}} \\ \hspace*{4cm} -
\frac{\left(|\T_{n}^{*}|^{-1} \sum_{i\in\T_{n}^{*}}
h_{1}(\Delta_{i})\right)\left(
|\T_{n}^{*}|^{-1}\sum_{i\in\T_{n}^{*}}
g_{2}(\Delta_{i})\right)}{B_{n}}$,

\vspace{0.25cm}
where $\displaystyle B_{n} =
|\T_{n}^{*}|^{-1}|\T_{n}^{1,0}|\left(|\T_{n}^{*}|^{-1}
\sum_{i\in\T_{n}^{*}} h_{2}(\Delta_{i})\right) -
\left(|\T_{n}^{*}|^{-1} \sum_{i\in\T_{n}^{*}}
h_{1}(\Delta_{i})\right)^{2}$.

Recalling (\ref{sample_sum}), we then have for all $\delta>0$ and
$a>0$
\begin{eqnarray}\label{di_gw}
\P\left(|\widehat{\alpha}_{0}^{n}-\alpha_{0}|>\delta|W\geq a\right)
\leq \P\left(\frac{|\T_{n}^{*}|^{-1}|\T_{n}^{1,0}|
|\overline{M}_{\T_{n}^{*}}(g_{1})|}{|B_{n}|}
> \frac{\delta}{2}\Big|W\geq a\right) \hspace{1cm} \\  +
\P\left(\frac{|\overline{M}_{\T_{n}^{*}}(h_{1})|
|\overline{M}_{\T_{n}^{*}}(g_{2})|}{|B_{n}|}
> \frac{\delta}{2}\Big|W\geq a\right). \notag
\end{eqnarray}
For the first term of the right hand of (\ref{di_gw}), since
$|\T_{n}^{*}|^{-1}|\T_{n}^{1,0}|\leq 1$, we have for all $\gamma>0$
\begin{eqnarray*}
\P\left(\frac{|\T_{n}^{*}|^{-1}|\T_{n}^{1,0}|
|\overline{M}_{\T_{n}^{*}}(g_{1})|}{|B_{n}|}
> \frac{\delta}{2}\Big|W\geq a\right) \leq
\P\left(|B_{n}|<\gamma|W\geq a\right) \hspace{2cm} \\ +
\P\left(|\overline{M}_{\T_{n}^{*}}(g_{1})|>\frac{\delta\gamma}{2}\Big|W\geq
a\right).
\end{eqnarray*}
Notice that

\vspace{0.25cm}

$\displaystyle B_{n} - (p_{1,0}^{2}\mu_{2} - p_{1,0}^{2}\mu_{1}^{2})
= p_{1,0}\mu_{2}\left(\frac{|\T_{n}^{1,0}|}{|\T_{n}^{*}|} -
p_{1,0}\right) +
\frac{|\T_{n}^{1,0}|}{|\T_{n}^{*}|}\overline{M}_{\T_{n}^{*}}(h_{2} -
p_{1,0}\mu_{2}) \\ \hspace*{4cm} -
\left(\overline{M}_{\T_{n}^{*}}(h_{1} - p_{1,0}\mu_{1})\right)^{2} -
2p_{1,0}\mu_{1}\overline{M}_{\T_{n}^{*}}(h_{1}-p_{1,0}\mu_{1})$

\vspace{0.25cm}

and

\vspace{0.25cm}

$\displaystyle \{|B_{n}|<\gamma\} \subset
\left\{|B_{n}-(p_{1,0}^{2}\mu_{2} - p_{1,0}^{2}\mu_{1}^{2})| >
|p_{1,0}^{2}\mu_{2} - p_{1,0}^{2}\mu_{1}^{2}| - \gamma\right\}$.

\vspace{0.25cm}

We then have for all $\displaystyle 0<\gamma <
\frac{2|p_{1,0}^{2}\mu_{2} - p_{1,0}^{2}\mu_{1}^{2}|}{2+\delta}$,

\vspace{0.25cm}

$\displaystyle \P\left(\frac{|\T_{n}^{*}|^{-1}|\T_{n}^{1,0}|
|\overline{M}_{\T_{n}^{*}}(g_{1})|}{|B_{n}|}
> \frac{\delta}{2}\Big|W\geq a\right) \\ \leq \P\left(|B_{n}|<\gamma|W\geq a\right)
+
\P\left(|\overline{M}_{\T_{n}^{*}}(g_{1})|>\frac{\delta\gamma}{2}\Big|W\geq
a\right) \\ \leq \P\left(|B_{n}-(p_{1,0}^{2}\mu_{2} -
p_{1,0}^{2}\mu_{1}^{2})| > \frac{\gamma\delta}{2}\Big|W\geq a\right)
+
\P\left(|\overline{M}_{\T_{n}^{*}}(g_{1})|>\frac{\delta\gamma}{2}\Big|W\geq
a\right) \\ \leq
\P\left(|\overline{M}_{\T_{n}^{*}}(g_{1})|>\frac{\delta\gamma}{2}\Big|W\geq
a\right) + \P\left(|\overline{M}_{\T_{n}^{*}}(h_{2}-p_{1,0}\mu_{2})|
> \frac{\gamma\delta}{8}\Big|W\geq a\right) \\ + \P\left(|\overline{M}_{\T_{n}^{*}}(h_{1}-p_{1,0}\mu_{1})|
> \frac{\sqrt{\gamma\delta}}{2\sqrt{2}}\Big|W\geq a\right) +
\P\left(\left|\frac{|\T_{n}^{1,0}|}{|\T_{n}^{*}|} - p_{1,0}\right| >
\frac{\delta\gamma}{8p_{1,0}\mu_{2}}\Big|W\geq a\right)
\\ + \P\left(|\overline{M}_{\T_{n}^{*}}(h_{1}-p_{1,0}\mu_{1})|
> \frac{\gamma\delta}{16p_{1,0}\mu_{1}}\Big|W\geq a\right)$.

\vspace{0.25cm}

From \cite{DY97}, Section 5, we have
\begin{equation}\label{eqd}
\P\left(\left|\frac{|\T_{n}^{1,0}|}{|\T_{n}^{*}|} - p_{1,0}\right| >
\frac{\delta\gamma}{8p_{1,0}\mu_{2}}\Big|W\geq a\right) =
\P\left(\left|\frac{1}{|\T_{n}^{*}|}\sum_{j=1}^{|\T_{n}^{*}|}
\left(T_{j} - p_{1,0}\right)\right| >
\frac{\delta\gamma}{8p_{1,0}\mu_{2}}\Big|W\geq a\right),
\end{equation}
where $\displaystyle \left(T_{j}\right)_{j\geq1}$ is a sequence of
i.i.d. Bernoulli random variables such that
\begin{equation*}
p_{1,0} = \P\left(T_{j} = 1\right) = 1 - \P\left(T_{j} = 0\right).
\end{equation*}

To majorize the right hand side of (\ref{eqd}), we use exactly the
same ideas that for the proof of Theorem \ref{thm:gw_main2} and Step
2 of the proof of Theorem \ref{thm:gw_main1}.

\medskip

For the second term of the right hand of (\ref{di_gw}), we have

\vspace{0.25cm}

$\displaystyle \P\left(\frac{|\overline{M}_{\T_{n}^{*}}(h_{1})|
|\overline{M}_{\T_{n}^{*}}(g_{1})|}{|B_{n}|}
> \frac{\delta}{2}\Big|W\geq a\right) \leq \P\left(
\frac{|\overline{M}_{\T_{n}^{*}}(g_{2})|}{|B_{n}|} >
\frac{\delta}{4p_{1,0}\mu_{1}}\Big|W\geq a\right) \\ \hspace*{1cm} +
\P\left( \frac{|\overline{M}_{\T_{n}^{*}}(g_{2})|}{|B_{n}|} >
\frac{\sqrt{\delta}}{2}\Big|W\geq a\right) +
\P\left(|\overline{M}_{\T_{n}^{*}}(h_{1} - p_{1,0}\mu_{1})| >
\frac{\sqrt{\delta}}{2}\Big|W\geq a\right)$

\vspace{0.25cm}

Now, the first and the second term of the right hand of the last
inequality can be treated as the first term of the right hand of
(\ref{di_gw}).

\medskip

Finally, to get the result, just apply Theorem \ref{thm:gw_main3} to
functions $g_{1}$, $g_{2}$, $h_{1}$ and $h_{2}$.

\vspace{15pt}

{\bf Acknowledgments.}  The author thanks its advisor Pr. Arnaud
GUILLIN for all his advices and suggestions during the preparation
of this work. \vspace{5pt}

\nocite{*}

\bibliographystyle{acm}
\bibliography{these}

\begin{thebibliography}{10}

\bibitem{AM}
{\sc Adamczak, R., and Mi{\l}o\'s, P.}
\newblock {CLT} for ornstein-uhlenbeck branching particle system.
\newblock {\em arXiv:1111.4559\/} (2011).

\bibitem{At94}
{\sc Athreya, K.}
\newblock Large deviations for branching processes-i. single type case.
\newblock {\em Ann. Appl. Probab. 5, No. 3\/} (1994), 779--790.

\bibitem{AtNe72}
{\sc Athreya, K., and Ney, P.}
\newblock {\em Branching Process}.
\newblock Springer, Berlin, 1972.

\bibitem{Azuma}
{\sc Azuma, K.}
\newblock Weighted sums of certain dependent random variables.
\newblock {\em Tôhoku Math. J 19, No. 3\/} (1967), 357--367.

\bibitem{Bennett}
{\sc Bennett, G.}
\newblock Probability inequalities for sum of independant random variables.
\newblock {\em Journal of the American Statistical Association 57, No. 297\/}
  (1962), 33--45.

\bibitem{BerSapGég}
{\sc Bercu, B., De~Saporta, B., and G\'egout-Petit, A.}
\newblock Asymptotic analysis for bifurcating autoregressive processes via a
  martingale approach.
\newblock {\em Electronic. J. Probab. 14\/} (2009), 2492--2526.

\bibitem{BiDj12}
{\sc Bitseki~Penda, S.~V., and Djellout, H.}
\newblock Deviation inequalities and moderate deviations for estimators of
  parameters in bifurcating autoregressive models.
\newblock {\em arXiv:1204.2355v1. Soumis\/} (2012).

\bibitem{BDG11}
{\sc Bitseki~Penda, S.~V., Djellout, H., and Guillin, A.}
\newblock Deviation inequalities, moderate deviations and some limit theorems
  for bifurcating {M}arkov chains with application.
\newblock {\em Preprint\/} (2011).

\bibitem{CS86}
{\sc Cowan, R., and Staudte, R.~G.}
\newblock The bifurcating autoregressive model in cell lineage studies.
\newblock {\em Biometrics 42\/} (1986), 769--783.

\bibitem{SGM12}
{\sc De~Saporta, B., G\'egout-Petit, A., and Marsalle, L.}
\newblock Asymmetry tests for bifurcating auto-regressive processes with
  missing data.
\newblock {\em preprint\/}, 2011.

\bibitem{DesaporGegMar}
{\sc De~Saporta, B., G\'egout-Petit, A., and Marsalle, L.}
\newblock Parameters estimation for asymmetric bifurcating autoregressive
  processes with missing data.
\newblock {\em Electronic Journal of Statistics 5\/} (2011), 1313--1353.

\bibitem{DelMar}
{\sc Delmas, J.~F., and Marsalle, L.}
\newblock Detection of cellular aging in {G}alton-{W}atson process.
\newblock {\em Stochastic Processes and their Applications 120\/} (2010),
  2495--2519.

\bibitem{DY97}
{\sc Dion, J., and Yanev, N.~M.}
\newblock Limit theorems and estimation theory for branching processes with
  increasing random number of ancestors.
\newblock {\em J. Appl. Probab. 34\/} (1997), 309--327.

\bibitem{Guyon}
{\sc Guyon, J.}
\newblock Limit theorems for bifurcating {M}arkov chains. application to the
  detection of cellular aging.
\newblock {\em Ann. Appl. Probab. Vol. 17, No. 5-6\/} (2007), 1538--1569.

\bibitem{Gu&Al}
{\sc Guyon, J., Bize, A., Paul, G., Stewart, E., Delmas, J., and Tadd\'ei, F.}
\newblock Statistical study of cellular aging.
\newblock {\em CEMRACS 2004 Proceedings, ESAIM Proceedings 14\/} (2005),
  100--114.

\bibitem{Hoeffding}
{\sc Hoeffding, W.}
\newblock Probability inequalities for sums of bounded random variables.
\newblock {\em Journal of the American Statistical Association 58, No. 301\/}
  (1963), 13--30.

\bibitem{SteMadPauTad}
{\sc Stewart, E.~J., Madden, R., Paul, G., and Tadd\'ei, F.}
\newblock Aging and death in an organism that reproduces by morphologically
  symmetric division.
\newblock {\em PLoS Biol 3(2)\/} (2005), e45.

\end{thebibliography}

\vspace{10pt}
\end{document}